\documentclass[ijds,fleqn,nonblindrev]{informs4}

\OneAndAHalfSpacedXI
\usepackage{natbib}
 \bibpunct[, ]{(}{)}{,}{a}{}{,}%
 \def\bibfont{\small}%
\makeatletter
\@ifundefined{setstretch}{
  
}{}
\makeatother
\usepackage{amsmath,amssymb} 
\usepackage{multirow} 
\usepackage{mwe}
\usepackage{esvect}
\usepackage{graphicx}
\usepackage[style=base]{caption}
\usepackage{stfloats} % for positioning of figure* on the same page
\usepackage{lipsum}
\usepackage{cite} 
\usepackage[caption=false]{subfig}
\usepackage{accents}
\usepackage{paralist, enumitem, nicefrac}
\usepackage{booktabs, tabularx, dcolumn, arydshln}
\usepackage{pgfplots,tikz,xcolor,adjustbox}
\usepackage[utf8]{inputenc}
\usepackage{pgf} 
\usepackage{arydshln}
\usepackage{color} 
\usepackage{todonotes} 
\usepackage{booktabs} %for table weighted rules 
\usepackage{nccmath}
\usepackage{float}
\usepackage{array}
\usepackage{mathrsfs}
\usepackage{blindtext}
\usepackage{algorithm}
\usepackage{algorithmic}
\newcolumntype{C}[1]{>{\centering\let\newline\\\arraybackslash\hspace{0pt}}m{#1}}
\newcommand{\ff}{\color{black}}

\TheoremsNumberedThrough     % Preferred (Theorem 1, Lemma 1, Theorem 2)
\ECRepeatTheorems

\EquationsNumberedThrough    % Default: (1), (2), ...
\MANUSCRIPTNO{IJDS}
\usepackage{booktabs} %for table weighted rules 
\begin{document}
%%%%%%%%%%%%%%%%

%\setcounter{page}{1}
%\DOI{10.1287/xxxx.0000.0000}%


\RUNTITLE{A Sensor-Driven Optimization Framework for Asset Management in Energy Systems}

% \TITLE{Bundling Information Goods of Decreasing Value}
\TITLE{A Sensor-Driven Optimization Framework for Asset Management in Energy Systems: Implications for Full and Partial Digital Transformation in Hydro Fleets}

\ARTICLEAUTHORS{%
%\AUTHOR{Snidely Slippery,\textsuperscript{a} Marg Arinella,\textsuperscript{b}}
%\AFF{\textsuperscript{a}Department of Bread Spread Engineering, Dairy University, Cowtown, IL 60208, \EMAIL{slippery@dairy.edu}; \textsuperscript{b}Institute for Food Adulteration, University of Food Plains, Food Plains, MN 55599, \EMAIL{m.arinella@adult.ufp.edu}} 
\AUTHOR{Farnaz Fallahi, Murat Yildirim}
\AFF{Industrial \& Systems Engineering, Wayne State University, Detroit, MI, USA, \EMAIL{farnaz.fallahi@wayne.edu, murat@wayne.edu}} %, \URL{}}
%\AUTHOR{Murat Yildirim}
%\AFF{Industrial \& Systems Engineering, Wayne State University, Detroit, MI, USA, \EMAIL{Murat@wayne.edu}}
\AUTHOR{Shijia Zhao, Feng Qiu}
\AFF{Argonne National Laboratory, Lemont, IL, USA, \EMAIL{zhaos@anl.gov,fqiu@anl.gov}}
%\AUTHOR{Feng Qiu}
%\AFF{Argonne National Laboratory, Lemont, IL, USA, %\EMAIL{fqiu@anl.gov}}
} % end of the block

\ABSTRACT{%
This paper proposes a novel prognostics-driven approach to optimize operations and maintenance (O\&M) decisions in hydropower systems. Our approach harnesses the insights from sensor data to accurately predict the remaining lifetime distribution of critical generation assets in hydropower systems, i.e., thrust bearings, and use these predictions to optimally schedule O\&M actions for a fleet of hydro generators. We consider complex interdependencies across hydro generator failure risks, reservoir, production, and demand management decisions. We propose a stochastic joint O\&M scheduling model to tackle the unique challenges of hydropower O\&M including the interdependency of generation capacities, the nonlinear nature of power production, operational requirements, and uncertainties. We develop a two-level decomposition-based solution algorithm to effectively handle large-scale cases. The algorithm incorporates a combination of Benders optimality cuts and integer cuts to solve the problem in an efficient manner. We design an experimental framework to evaluate the proposed prognostics-driven O\&M scheduling framework, using real-world condition monitoring data from hydropower systems, historical market prices, and water inflow data. The developed framework can be partially implemented for a phased-in approach. Our experiments demonstrate the significant benefits of the sensor-driven O\&M framework in improving reliability, availability, effective usage of resources, and system profitability, especially when gradually shifting from traditional time-based maintenance policies to condition-based prognostics-driven maintenance policies. 

%\footnote{ This research employs various data sources, including real-world degradation data for bearing systems and historical data for water inflows, electricity demand, and market prices. The authors recognize the importance of ethical considerations in data collection and use. All data sources used in this study are either publicly available or obtained with proper permissions and in compliance with relevant regulations and guidelines.}
% Current word count: 277 words
}%

\FUNDING{This material is based upon work supported by the U.S. Department of Energy, Office of Energy Efficiency and Renewable Energy, specifically the Water Power Technologies Office, titled “An Integrated Framework for Condition Monitoring based Asset Management in Hydropower Fleets:  Adaptive and Scalable Prognostics, Operations, and Maintenance''}

%\KEYWORDS{deterministic inventory theory; infinite linear programming duality;
%  existence of optimal policies; semi-Markov decision process; cyclic schedule}

 %\HISTORY{This paper wasfirst submitted on April 12, 1922 and has been with the authors for 83 years for 65 revisions.}

\maketitle
%%%%%%%%%%%%%%%%%%%%%%%%%%%%%%%%%%%%%%%%%%%%%%%%%%%%%%%%%%%%%%%%%%%%%%

% Samples of sectioning (and labeling) in MNSC
% NOTE: (1) \section and \subsection do NOT end with a period
%       (2) \subsubsection and lower need end punctuation
%       (3) capitalization is as shown (title style).
%
%\section{Introduction.}\label{intro} %%1.
%\subsection{Duality and the Classical EOQ Problem.}\label{class-EOQ} %% 1.1.
%\subsection{Outline.}\label{outline1} %% 1.2.
%\subsubsection{Cyclic Schedules for the General Deterministic SMDP.}
%  \label{cyclic-schedules} %% 1.2.1
%\section{Problem Description.}\label{problemdescription} %% 2.

% Text of your paper here

% Text of your paper here
\section{Introduction}
%with the degradation of {\ff generation assets} and prevent unexpected failures \citep{de2016review}. 
%Evidently, O\&M is a significant cost driver for energy systems, 
%In industrial systems, especially power systems, the management of asset maintenance is a significant cost driver \citep{de2020review}. Maintenance operations
%e.g. accounting for 20\%-25\% of the total levelized cost per kWh of wind turbines \citep{el2012operation}.
%To enable a succesful, a simple focus on reducing maintenance costs or improving reliability is insufficient. 
Operations and maintenance (O\&M) is a fundamental problem in energy systems, with widespread implications for cost and reliability \citep{froger2016maintenance}. A central objective for O\&M is to ensure operational availability and reliability of energy assets by performing maintenance actions as needed \citep{de2016review}. %in an effort to %Achieving this objective requires a delicate 
%balance maintenance downtime, production schedules, and the delivery of cost-effective and reliable power \citep{ding2015maintenance}. 
%O\&M considers both operations (i.e. production) and maintenance scheduling. 
One of the key trade-offs in maintenance scheduling revolves around determining the frequency of maintenance actions. Opting for frequent maintenance intervals can reduce equipment's useful life and disrupt operations due to maintenance outages; while delaying maintenance increases the risk of failures, which can have severe repercussions for grid stability. Consequently, operators are constantly seeking methods to enhance their awareness of asset conditions and failure risks to execute maintenance promptly. Leveraging sensor-driven insights facilitated by digital transformation has proven instrumental in this regard, enabling operators to monitor real-time changes in asset conditions and predict failures. However, when operations and maintenance are integrated into a joint model, the trade-offs become more intricate and multifaceted, complicating the translation of failure predictions into actionable prescriptive insights. This complexity is particularly pronounced in hydropower systems, where operational interactions, such as reservoir management, production curves, and energy prices, introduce numerous interdependencies across assets and decision-making processes. This disparity between prediction and prescription, often called the  ``predictive-prescriptive gap", represents a significant obstacle to fully realizing the benefits of digital transformation for O\&M outcomes. There is a need for a new generation of O\&M models that can intergrate sensor-driven insights on asset conditions within large-scale O\&M models capable of modeling fleet-level interactions.%from a fleet of assets within large-scale operational models that explicitly model the complexity of operational interdependencies and consequences. 

\indent
Hydropower industry - a significant source of global renewable energy production - is no exception to this need \citep{keizer2017condition,/content/publication/2ef8cebc-en}. 
Maintenance breakdowns of {\ff hydro generation assets} have exceeded acceptable limits, causing concern among stakeholders \citep{efficiency2016water}. Factors such as fluctuating production levels to satisfy operational and environmental requirements contribute to the degradation of critical {\ff assets such as generators and turbines} and reduce the intervals between maintenance activities \citep{yang2016wear,saarinen2017allocation,liu2016review}. Frequent outages of {\ff hydro generation assets} lead to forced spillage of water with zero or reduced economic value, generation loss, unsatisfied demand, and, in some cases, financial penalties.  These circumstances underline the importance of proper maintenance management policies and preventive actions to ensure high system availability and reliability.
% Excessive vibration, premature failure, and damage to critical components of these degraded turbines have been reported and have led to the shutdown of some plants for months \citep{liu2016review}. Proper maintenance management policies and preventive actions can potentially prevent failures and ensure high availability and reliability of the system.
Developing cost-efficient maintenance plans for {\ff hydro generation assets} requires consideration of several factors.\citep{helseth2018optimal}. These include the operational interdependency and the varying generation capacity of {\ff hydropower systems}, which depend on {\ff generation assets} availability, natural water inflow, and reservoir storage levels \citep{helseth2018optimal, labadie2004optimal}. Consequently, effective maintenance schedules must be proactively planned, remain adaptive to changing system conditions and operational requirements, and flexible to respond to uncertainties.

% Another factor is the varying generation capacity of hydropower plants, which depends on the random natural water inflow and the amount of stored water in reservoirs. The connection of reservoirs along a river basin could impose additional complexities by the dependencies among 
% their generation capacities \citep{labadie2004optimal}. It is vital to predict and evaluate the effect of outages of major units on upstream or downstream hydropower plants in connected reservoirs. Because of these factors, effective maintenance schedules must be planned to ensure that the hydroturbines run safely and reliably within their design lifetime. More importantly, to guarantee the successful implementation of the scheduled O\&M activities, these decisions should evolve in response to the continuous changes in the system conditions and operational requirements and should provide enough flexibility to react to the underlying uncertainties.

\indent
Various maintenance management policies have been proposed to tackle O\&M challenges in energy systems. One of the widely-adopted maintenance approaches, called \textit{reactive} maintenance policy, relies on corrective maintenance actions following {\ff generation asset} breakdowns. In an effort to reduce failures, especially for critical {\ff generation assets}, \textit{preventive} maintenance policies are developed, where the maintenance on {\ff assets} is performed following predetermined intervals. A widely used preventive maintenance strategy in the energy industry, especially for {\ff hydropower systems}, is the time-based maintenance (TBM) policy, which schedules maintenance based on calendar time or machine run-time intervals \citep{guedes2015continuous,osti_1330494}. While this maintenance policy offers some advantages, it often results in unnecessary maintenance actions, high maintenance costs, and does not prevent catastrophic failures. These maintenance policies are developed based on the assumption that the degradation, and consequently failure, behavior of all {\ff generation assets} of the same type are similar. However, each individual {\ff asset}, to some extent, possesses unique degradation, tear and wear trends. Limited flexibility and strong assumptions on failure risks, reduces the effectiveness of TBM policies for achieving the reliability and availability requirements.
\newline
\indent
The evolution of condition monitoring (CM) technologies, such as sensors and computational devices, has significantly enhanced the reliability and efficiency of maintenance policies. These technologies utilize an array of sensors to continuously or periodically record and monitor the health indicator parameters of {\ff individual generation assets} during operation. The information gathered by these sensors is transformed into unique identifiers known as \textit{degradation signal}s, which are then used to schedule \textit{condition-based} maintenance (CBM) actions. Implementation of CBM policies, as indicated by independent surveys across various industries, can potentially reduce maintenance costs by up to 30\% and equipment downtime by as much as 45\% \citep{sullivan2010operations}. CBM policies rely on CM-based failure predictions, which take one of two forms: diagnostics and prognostics. Diagnostic approaches use sensor data to estimate the current state of health \citep{kusiak2011prediction,hameed2009condition}, and are typically used to identify {\ff assets} with imminent failure risks. Prognostic approaches derive remaining life distributions (RLDs) for {\ff generation assets}, which require an estimate of the current state of health (as in diagnostics), as well as an accurate prediction of how the health state of the {\ff asset} is likely to evolve in the future \citep{selak2014condition,wang2016review}. 

In line with the existing CM approaches, the CBM approaches in power systems also focus on two types of policies: diagnostics-driven, and prognostics-driven maintenance policies.  In diagnostics-driven policies, CM systems alert the operators when the failure risks of {\ff generation assets} reach a certain severity. Depending on the subjective judgment of the maintenance personnel, these alerts are often used to initiate immediate maintenance actions. Such diagnostics-driven CBM policies rely on the current state of {\ff the assets} and do not provide advanced notice for planning maintenance actions. To overcome this constraint, real-time sensor information can be synergized with analytical methods to accurately predict the progression of faults and the remaining life of the {\ff assets}. If the degradation trends of {\ff the asset} demonstrate predictable patterns and can be effectively modeled, they can be used to predict the remaining life distribution (RLD) of the {\ff asset}. This information paves the way for the development of a prognostics-driven CBM approaches. Prognostics-driven CBM approaches introduce the flexibility and agility necessary to create a well-balanced O\&M schedule that fulfills reliability and availability requirements while also meeting the operational objectives of the system.

Despite the promising benefits of sensor-driven CBM maintenance policies and supporting technologies, their practical implementation often faces hurdles. High initial investment costs and long payback periods can deter industries from transitioning to CBM policies. Even when the required technologies are implemented, the full potential of collected sensor data may not be realized due to the complexities in decision integration. In the face of uncertain operational environments, aligning maintenance and operational decisions is also challenging and hasty responses can compromise both operational and maintenance objectives. Maintenance policies should therefore allow sufficient time to formulate the most effective response. This necessitates early predictions of {\ff generation assets} failures to enable proactive decision-making. However, existing CM technologies are typically focused on monitoring and diagnostics, falling short in predicting failure risks. Moreover, it's not sufficient to merely make asset-specific failure risk predictions and maintenance decisions; a comprehensive understanding of the entire system's performance is essential. Take the example of a situation in {\ff hydropower systems} where a key upstream {\ff generation asset} is under maintenance, potentially limiting available water flow to other downstream assets. In this case, there may be additional incentives to conduct maintenance on a downstream assets at the same time, to better utilize fleet-level operational capacity. 
%, the impact of also taking a downstream {\ff generation asset} offline might be lessened. 
Therefore, decision-making should extend beyond individual {\ff asset} O\&M decisions; making a comprehensive fleet-level performance overview indispensable. This holistic perspective facilitates fleet-level alignment, leading to substantial cost reductions.  Without thoughtful integration of failure predictions into decision-making tools, there is a risk of poor maintenance decisions, leading to discouragement and policy disuse. These challenges underscore the need for more robust implementation of prognostics-driven CBM policies.
\indent
 This study explores the benefits of prognostics-driven CBM
 for joint O\&M planning in {\ff hydropower systems}. {\ff The hydropower systems} consist of multiple {\ff hydropower plants} with various generation and reservoir capacities.  Our novel sensor-driven framework enhances the understanding of utilizing real-time CM information to manage failure risks of {\ff generation assets, i.e.,  hydro generators,} and coordinate the backbone functions of hydropower systems, specifically maintenance and operations activities. The main contributions of this paper are fourfold.
 \begin{itemize}
 \item We develop a novel prognostics-driven framework for O\&M scheduling of {\ff hydropower systems}. The framework uses sensor-driven prognostics to drive optimal O\&M schedules for {\ff a fleet of hydro generators}. The proposed framework captures the dynamic interactions between (i) the failure risks and remaining lifetime predictions of {\ff hydro generators}, and (ii) operations of {\ff hydropower fleets}. The proposed model ensures that the full potential of sensor-driven prognostics can be unlocked to increase system-level efficiency and reliability in hydro fleets. 
\item We construct a stochastic joint O\&M scheduling model specifically for {\ff hydropower systems}. The developed model addresses some of the unique characteristics of {\ff hydropower system} O\&M including the spatial and temporal interconnections between reservoirs, the nonlinearity of power production, and the underlying operational uncertainties, requirements, and system interdependencies. Uncertainties associated with water inflows and market price uncertainties are considered in the O\&M scheduling through the notion of scenarios.
\item We propose a two-level decomposition-based solution algorithm to solve the mixed-integer stochastic O\&M scheduling model for large-scale cases. We exploit the structure of the model to reformulate the problem in a two-stage formulation. Leveraging on the decomposed model, we propose an algorithm that uses a combination of Benders optimality cuts and integer cuts, to obtain optimal solutions in an efficient manner.
\item We design an experimental framework that allows to evaluate the performance of hydropower systems based on the common key maintenance, operations, and cost metrics within the hydro industry. The developed framework builds on real-world hydropower sensor data acquired from an industrial collaborative called Hydropower Research Institute (HRI) to mimic the degradation of {\ff hydro generators} and drive remaining life estimations using Bayesian statistics. We utilize historical market price and water inflow data to represent the operational uncertainties of hydropower systems. The framework allows partial implementation of prognostics-driven CBM maintenance policy, i.e., managing a subset of {\ff hydro generators} under the CBM policy. 
\end{itemize}
\indent
We conduct a series of in-depth experiments on {\ff connected hydropower plants}, depicted in Fig \ref{fig:Casestudy}, to analyze the performance of the sensor-driven CBM framework in a realistic setting. Our results reveal the significant merits of the sensor-driven framework in improving several aspects of {\ff hydropower systems}, such as reliability, availability, effective use of resources, and system profitability, and we conduct a unique set of experiments to investigate the impact of a gradual shift from the conventional TBM policy, which does not consider real-time sensor information, to the CBM policy proposed in this paper.
\newline
\indent
The remainder of this paper is organized as follows. The subsequent section provides a review of the relevant literature. In Section \ref{sec:sensorframework}, we detail our sensor-driven prognostics CBM framework, covering the degradation process of {\ff hydro generators}, the continuous predictions of their RLDs, and its connection to the maintenance cost functions. We then introduce our Prognostics-Driven Operations and Maintenance Scheduling (POMS) model in Section \ref{sec:model}, offering an in-depth exploration of operational and maintenance constraints inherent to the architecture of {\ff hydropower systems}. Section \ref{sec:Algorithm} proposes a strategic reformulation of the POMS model and presents a decomposition-based solution algorithm designed to address the  computational scalability of large-scale instances. The experimental results and key findings are presented in Section \ref{sec:experiments}. The paper concludes with reflections on our research, outlining potential directions for future work in the final section.
\section{Literature}
In this section, we provide a brief overview of the related literature on condition monitoring and prognostic approaches, the CBM literature on the O\&M planning in power systems, and the relevant studies on the O\&M scheduling problem with a specific focus on hydropower systems.
\newline
\indent
\subsection{Condition Monitoring and Predictive Analytics}
Integrated CM systems are actively used for monitoring the performance
of various {\ff generation assets}, by collecting and analyzing data
from sensors. Proximity and eddy current probes as well as vibration and temperature sensors are among the most widely used sensors.  
\begin{figure}[htbp]
    \centering
    \includegraphics[scale=0.4]{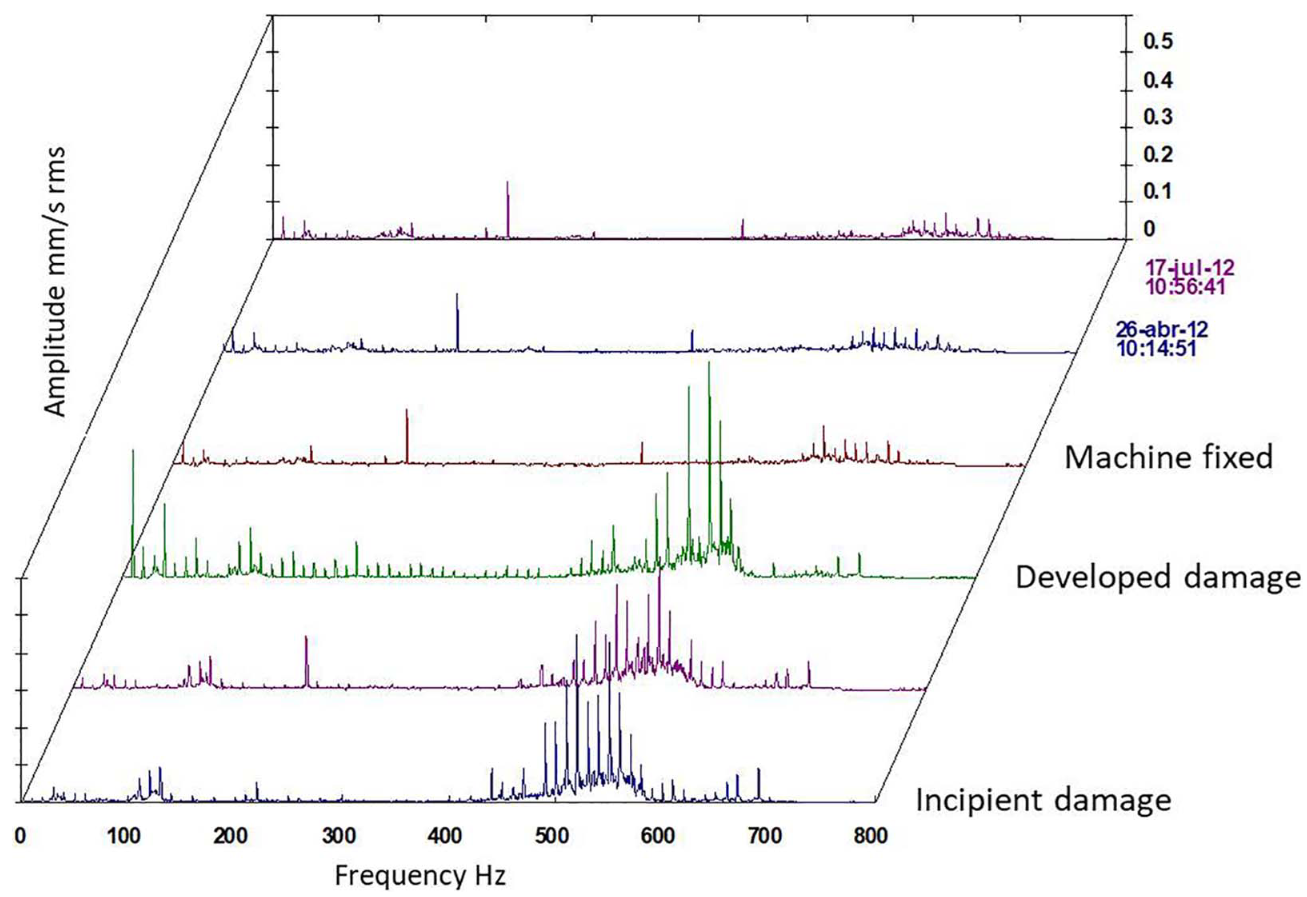}    \caption{Evolution of the vibration signal of a hydro generator before and after maintenance \citep{egusquiza2018advanced}}
    \label{fig:Vibration}
\end{figure} 
In {\ff hydropower systems}, operators generally use sensor data for detecting anomalies and tripping units for system protection \citep{kougias2019analysis}.  The focus is mainly on monitoring major equipment, such as turbines and generators \citep{efficiency2020water,ribeiro2014equipment}. Fig \ref{fig:Vibration} illustrates the progression in the vibration signal of {\ff a hydro generator} before and after maintenance. The high amplitude of vibration, 500-600 Hz, indicates the possibility of incipient damage \citep{standard1996mechanical}. Existing CM techniques can be classified into mechanical, electrical, and chemical or hydraulic monitoring,
depending on the specific application \citep{tavner1987condition,klempner2004operation}. Mechanical monitoring includes vibration analysis, shaft displacement
measurements, and bearing temperature monitoring, which can help detect issues such as
misalignment, imbalance, and bearing wear. Electrical monitoring focuses on the
insulation system, including partial discharge measurements, and the analysis of electrical
parameters like current, voltage, and power factor to identify winding issues, rotor faults,
and stator core faults. Hydraulic monitoring involves observing parameters such as pressure, flow, and temperature within the {\ff hydro generators} to detect issues related to
the turbine, such as cavitation, flow instabilities, and wear of turbine components. For reviews of the advancements in fault diagnostics and prognostics approach for {\ff hydro
equipment}, readers are referred to \citep{wang2016review,selak2014condition}

\subsection{CBM for Operations and Maintenance Planning in Power Systems} 
The majority of prognostics-driven CBM approaches in energy systems focus on {\ff single-unit systems}\citep{de2020review}. {\ff A single-unit system is composed of a single asset or unit}. Markov chain models are used to characterize degradation and derive optimal maintenance decisions \citep{mahani2018joint}. More examples of CBM for single-unit systems can be found in \citep{besnard2010approach,byon2010optimal,byon2010season}. This literature is expanded for more complex settings by modeling the system structure and introducing various types of interdependencies between the {\ff generation assets}, such as economic, stochastic, structural, and resource dependencies\citep{keizer2017condition,shin2015condition,alaswad2017review}. Recently, \citep{yildirim2016sensor} proposed a fleet-wide CBM policy for generator maintenance scheduling in vertically integrated power systems. While benefiting from prognostic predictions, these models do not capture the interdependency between the O\&M actions. Scheduling maintenance activities in isolation from operations and/or actual characteristics of the system results in CBM policies that are sub-optimal, not reliable, and efficient enough.  Any viable CBM policy needs to consider {\ff the interdependency within the system} when scheduling maintenance
actions to reduce production and economic loss \citep{helseth2018optimal}.

In comparison to those on pure maintenance scheduling problems, the number of studies on joint O\&M scheduling in power systems is limited. Some of the researchers investigated  the joint O\&M scheduling problems in wind farms \citep{yildirim2017integrated,bakir2021integrated,papadopoulos2023joint}, vertically integrated grid \citep{basciftci2018stochastic}, and microgrids \citep{fallahi2021predictive}. The problems are modeled either by deterministic or stochastic approaches. Typical approaches to connecting prognostics information to O\&M decisions include maintenance cost functions \citep{yildirim2016sensorII}, uncertainty sets \citep{basciftci2018stochastic}, and chance constraints \citep{basciftci2020data}. These studies show the effectiveness of sensor-driven O\&M scheduling approaches in improving the reliability, revenue, and other operational metrics of {\ff the energy system} under consideration \citep{bakir2021integrated,yildirim2019leveraging}.
\subsection{Operations and Maintenance Planning in Hydropower Systems}
Many studies in the literature have addressed the main challenges in hydropower operations scheduling \citep{singh2017operation,labadie2004optimal,steeger2014optimal}.
The operations scheduling of {\ff hydropower plants} are characterized as a complex problem due to its non-linearity, non-convexity, stochasticity, and, in most cases, large-scale.
% The non-linearity and non-convexity arise mainly from the relationship between the power output of the unit with the turbine water discharge, performance efficiencies, and the level of stored water in the feeding reservoir (water head effect) \cite{chang2001experiences}.
To deal with the non-linearity of power production, several methods have been proposed, including piecewise linear approximations \citep{conejo2002self,sarasty2018mixed} and smoothing splines \citep{seguin2015self}. Some studies explore modeling simplification techniques such as dropping the power production dependency on the net water head when reservoirs or the length of the planning horizon is large enough \citep{dashti2016weekly}. Another modeling technique is the plant-based power production concept in which the power production is modeled in an aggregate \citep{rodriguez2018milp}, i.e., plant-based, manner rather than the individual \citep{conejo2002self}, i.e., {\ff asset-based}, manner. Although this method reduces the problem size and decreases the computational time \citep{kong2020overview}, it ignores the precise representation of individual {\ff hydro generator} characteristics. Accurate O\&M planning of {\ff hydropower systems}, however, depends highly on the individual {\ff assets} operational characteristics and their interactions with each other. Some studies also give great importance to the necessity of counting the start-up costs of the {\ff hydro generators} to reduce water loss and equipment wear \citep{nilsson1997hydro,dashti2016weekly} which require precise {\ff asset-wise} on/off status consideration. Modeling the underlying uncertainties within {\ff hydropower systems} such as energy generation capacity and market conditions is also vital for planning O\&M activities efficiently and in a reliable manner. In addition to the electricity prices volatility, the natural water inflows into the reservoirs are another main source of uncertainty in the {\ff hydropower systems}. Two-stage stochastic programming \citep{pousinho2012short,baslis2011mid,khodayar2013enhancing,egging2016linear}, multi-stage stochastic programming, chance-constrained programming \citep{van2014joint}, robust optimization \citep{dashti2016weekly}, and dynamic programming \citep{arce2002optimal} are some of the modeling approaches proposed in the literature to address these operational uncertainties. 

In comparison to the literature on operations scheduling, little attention has been paid to the joint O\&M planning in {\ff hydropower systems}. The existing studies focus more on the operational side, and maintenance activities are generally managed by enforcing additional constraints\citep{guedes2015continuous,sarasty2018mixed,helseth2018optimal,ge2018mid}. With the advancement of CM technologies, monitoring {\ff hydro generators} become one of the best practices in managing {\ff assets} health condition \citep{efficiency2020water,ribeiro2014equipment}. There are very few studies that employ the CBM policy, i.e., scheduling maintenance actions based on sensor data and the {\ff asset'} state of health. Most studies utilize diagnostics-driven CBM approaches and consider a {\ff single-asset hydropower system} \citep{fu2004predictive,betti2021condition}. The focus of these studies is on fault detection ability, and they usually have no or a very simplified maintenance decision-making optimization approach.  Very few studies focus on prognostics-driven CBM approaches to managing maintenance actions. \citep{qian2014condition} proposes a CBM policy for a {multi-component hydro generator} and uses a proportional hazard model to represent the risk of failure. Some studies have used the Markov decision process to model the CBM optimization problem in {\ff hydropower systems} \citep{welte2008deterioration}.
\newline
\indent
One of the common key limitations of these studies is that the main maintenance and/or operational characteristics of {\ff hydropower systems} are either ignored or oversimplified. 
To date, the proposed research on joint O\&M scheduling in {\ff hydropower systems} has not captured and investigated the potential of CBM policies in the complex operations of {\ff hydropower systems}. 

\section{ Sensor-driven Degradation Modeling and Prognostics Framework}\label{sec:sensorframework}
In this section, we introduce the degradation modeling and sensor-driven failure prognostics framework for {\ff hydro generators}. The developed framework is composed of four main steps. The first step models the degradation signals of {\ff a population of hydro generators} using a parameterized stochastic model. 
The second step uses a Bayesian updating procedure and leverages the sensor data to improve the accuracy of the initial degradation modeling attempts. More specifically, to capture the {\ff unit-to-unit} degradation variability and improve the accuracy of degradation modeling, the Bayesian procedure combines the {\ff population-level knowledge} of the {\ff hydro generators} degradation process with the real-time sensor data collected from individual assets while operating. Using the updated degradation models, in the next step, the failure predictions of {\ff individual hydro generators} are reevaluated. 
Finally, the last step utilizes the most recent prediction of RLD to calculate a prognostics-driven maintenance cost function of {\ff hydro generators} over time. 

Before modeling the degradation signals, we first elucidate some important terminologies. Let $X_i(t)$ denote the cumulative degradation of {\ff hydro generator} $i$ up to time $t\in\mathbb{R}_+$. During the operation of {\ff hydro generators}, degradation signals are collected and processed by the CM system at discrete times $t_i^0,t_i^1,...,t_i^o$. We represent the observed values by $x_i^0,x_i^1,...,x_i^o$ where $x_i^k=X_i(t_i^k)$.
We define failure occurrence as the first time when the amplitude of the accumulated degradation exceeds a predetermined failure threshold $\Lambda_i$, i.e., $\tau_i=\inf\{t\geq0|X_i(t)>\Lambda_i\}$. Consequently, the remaining life of the {\ff hydro generator} $i$ is the minimal time it takes until the degradation signal crosses the failure threshold. Given the hydro-generator survival until the observation time $t_i^o$, the remaining life of the {\ff asset} is mathematically defined as follows:
\begin{equation}
  R_i^{t_i^o}=inf\{s\geq0|X_i(s+t_i^o)\geq\Lambda_i|X_i(t)<\Lambda_i, \forall t\leq t_i^o\}
  \label{eq:RLV}
\end{equation}
We now proceed with introducing our degradation modeling approach and using sensor observations to improve RLD estimations of individual {\ff hydro generators}.
\subsection{Prognostics - Degradation Modeling and Remaining Life Predictions}
\label{sec:DMF}
There are different ways to describe the trajectory of  {\ff hydro generators} degradation over time \citep{zhang2018degradation}. In this paper, we model the degradation as a continuous-time continuous-state stochastic process as follows:
\begin{equation}
  X_i(t)=\eta_i(t;\phi,\beta_i)+\epsilon_i(t;\sigma)
  \label{eq:Deg}
\end{equation}
\indent
In this formulation, $\eta_i(.)$ represents the general functional form of the degradation signal along with $\phi$ and $\beta_i$ as two independent explanatory parameters. The parameter $\phi$ describes the degradation characteristics common between a {\ff population of hydro generators}, i.e.,  generation assets of the same type, and is generally deterministic and constant. The degradation variability from {\ff asset to asset} is captured by parameter $\beta_i$ which is stochastic and follows a statistical distribution across the population of {\ff hydro generators}. The error term $\epsilon_i(.)$ represents the noise level due to degradation and measurement errors through the variance parameter $\sigma$. 

In practice, for modeling particular degradation processes, experts use domain knowledge and historical data to estimate the value of the deterministic parameters and the distribution of the stochastic parameters $\beta_i$ and $\epsilon_i$. The resulting population-based degradation model can be updated to specialized ones by incorporating the improved estimations of parameters that are degradation-path dependent,{\ff i.e., $\beta_i$}. In other words, we use the degradation indicator sensor data collected from {\ff assets} to periodically reevaluate the stochastic degradation model parameters. To begin with, we denote the initial estimation of the stochastic parameter by $\pi(\beta_i)$ and consider a sequence of observed degradation signals by the time of observation $t_i^o$. Given the partial degradation signal observations $x_i^0,x_i^1,...,x_i^{t_i^o}$, we obtain the conditional distribution of the stochastic parameter through the following Bayesian framework. 
\begin{equation}
P(\beta_i|x_i^0,x_i^1,...,x_i^{t_i^o})=\frac{P(x_i^0,x_i^1,...,x_i^{t_i^o}|\beta_i)\pi(\beta_i)}{P(x_i^0,x_i^1,...,x_i^{t_i^o})}
  \label{eq:Bayes}
\end{equation}
where $ P(\beta_i|x_i^0,x_i^1,...,x_i^{t_i^o})$ is the posterior distribution of the stochastic parameter $\beta_i$ at the observation time, and the denominator is a constant normalization factor. We denote this posterior distribution of the stochastic parameter by $u(\beta_i|x_i^0,x_i^1,...,x_i^{t_i^o})$. Once the distribution of the stochastic parameter and the updated degradation models are acquired, we reevaluate the remaining life of the {\ff hydro generators}. 
The conditional distribution of the remaining life of the hydrogenerator $i$, $R_i^{t_i^o}$, given the observed values $x_i^0,x_i^1,...,x_i^{t_i^o}$ are as follows:
\begin{align}
  F_{R_i^{t_i^o}|x_i^0,x_i^1,...,x_i^{t_i^o}}(t)&=P(R_i^{t_i^o}\leq t|x_i^0,x_i^1,...,x_i^{t_i^o})=P\{X_i(t+t_i^o)\geq\Lambda_i|x_i^0,x_i^1,...,x_i^{t_i^o}\}
  \label{eq:RLD}
\end{align}
\indent
In some cases, a closed-form solution to the RLD estimation \eqref{eq:RLD} may exist \citep{gebraeel2005residual}. The frequency of the Bayesian updating procedure could vary among different industrial applications based on the cost of data acquisition.
\begin{figure}[!htbp]
    \centering
    \includegraphics[scale=0.4]{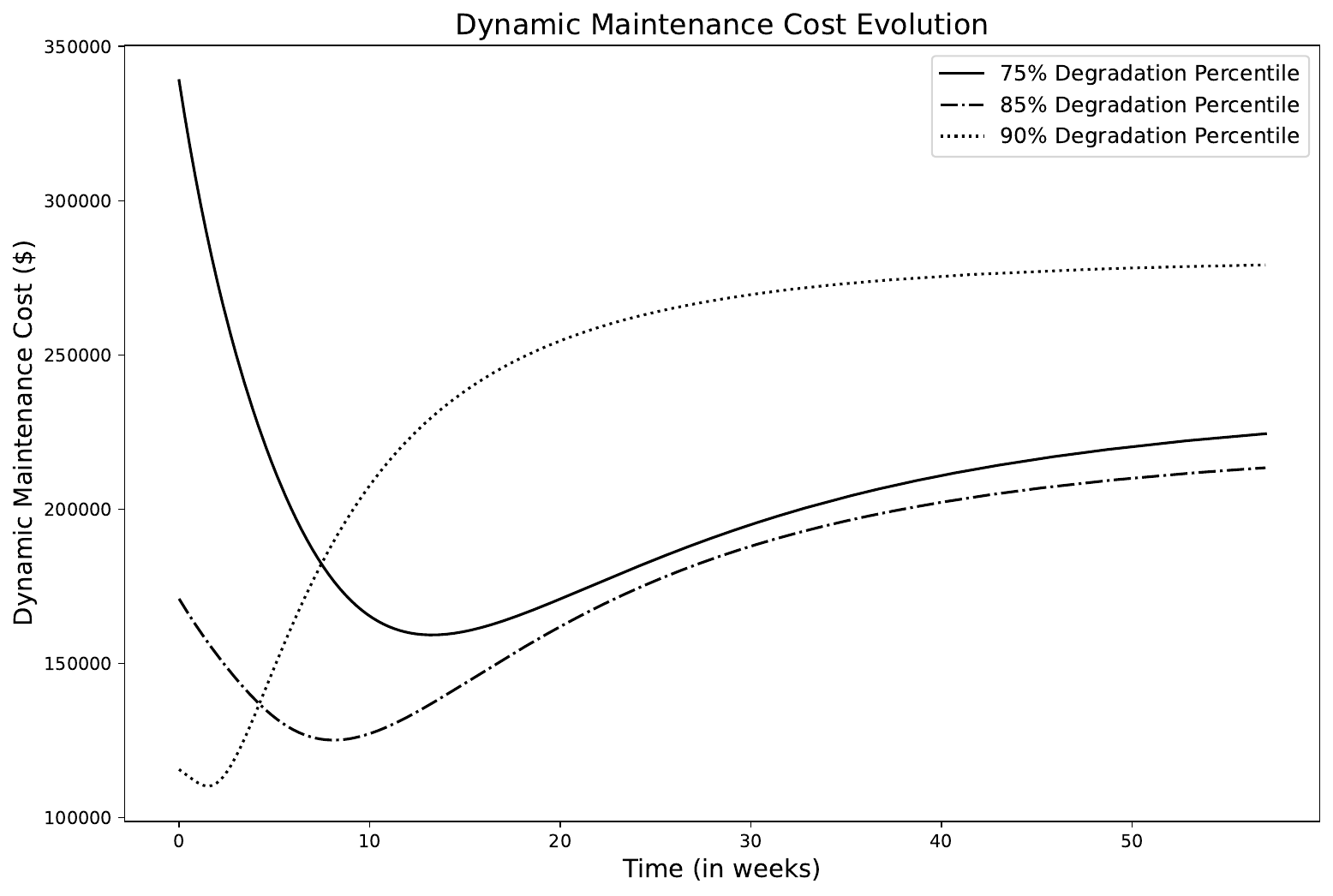}
    \caption{Degradation-based maintenance cost function of a hydro generator}
    \label{fig:costfunction}
\end{figure}  
\subsection{Degradation-based Maintenance Cost Function}\label{sec:DMC}
To link the sensor-driven RLD estimations {\ff of hydro generators} with the O\&M decisions of hydropower systems, we develop a maintenance cost function using renewal reward processes. Given the updated RLD of {\ff hydro generators}, we calculate the expected cost of conducting maintenance $t$ time periods after the observation time $t_i^o$ as follows \citep{yildirim2016sensor}: 
\begin{equation}
  C^{d,t_i^o}_{i,t} = \frac{ C^{p}_{i} P(R^{t_i^o}_{i} > t) + C^{f}_{i} P(R^{t_i^o}_{i} \leq t) }{\int_0^{t} P(R^{t_i^o}_{i} > z) dz + t_i^o}
  \label{eq:C_d}
\end{equation}
where $C^{p}_{i}$ and $C^{f}_{i}$ denote the maintenance costs associated with preventive and corrective maintenance actions, respectively. This cost function balances the trade-off between the risk of unexpected failures of {\ff hydro generators} and cost of preventive maintenance actions by combining the expected preventive maintenance costs per cycle $C^{p}_{i} P(R^{t_i^o}_{i} > t)$, the expected corrective maintenance costs per cycle $C^{f}_{i} P(R^{t_i^o}_{i}\leq t)$, and the expected maintenance cycle length $(\int_0^{t} P(R^{t_i^o}_{i} > z) dz + t_i^o)$ of the {\ff hydro generator}. The introduced maintenance cost functions allow us to integrate the most accurate RLD predictions {\ff of hydro generators} into the optimization model and ensure that the O\&M decisions of {\ff hydropower systems} adapt to failure risks and associated uncertainties of these {\ff assets}.
     \begin{figure}[!htbp]
    \centering
    \includegraphics[scale=0.5]{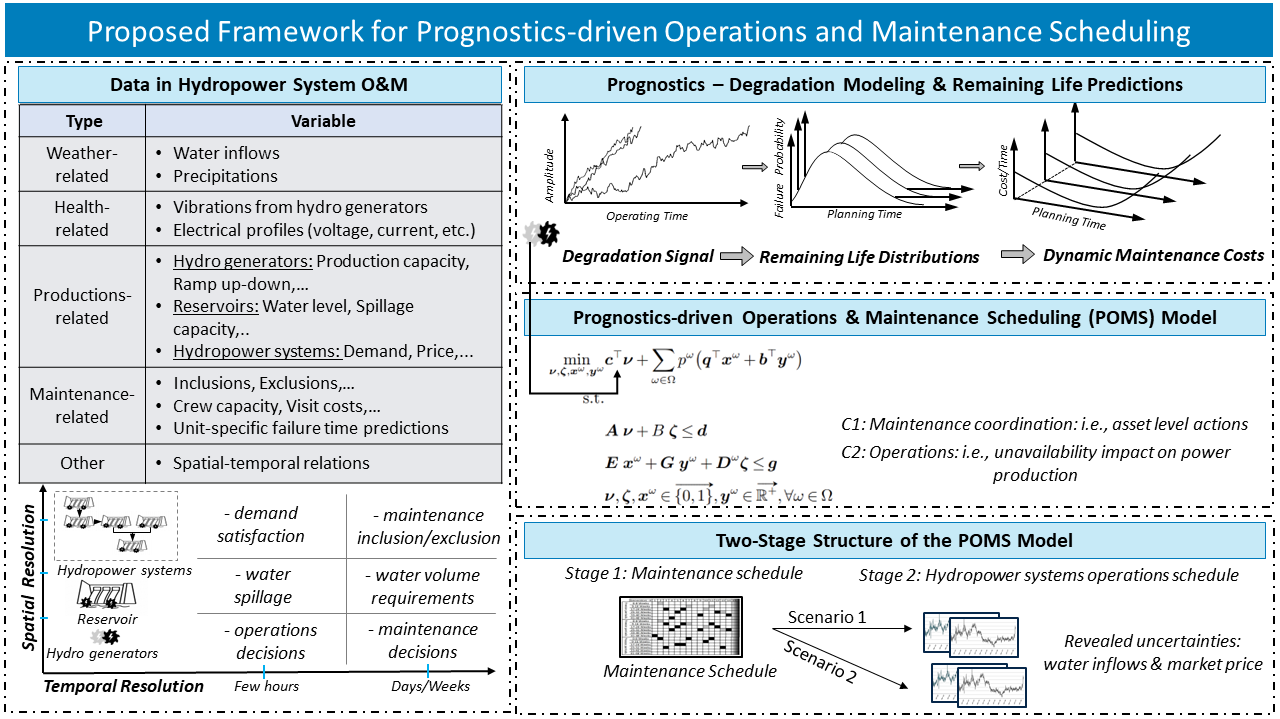}
    \caption{Prognostics-driven operations and maintenance scheduling framework} 
    \label{fig:Framework}
\end{figure}  

\section{Prognostics-driven CBM model for Joint O\&M Scheduling}\label{sec:model}
This section introduces a prognostics-driven Operations and Maintenance (O\&M) scheduling model for {\ff hydropower systems}, herein referred to as POMS. The problem is modeled as a stochastic mixed-integer linear programming problem, with the goal of optimizing maintenance actions for {\ff hydro generators} within the physical and operational constraints of {\ff hydropower systems} based on sensor-driven failure ris predictions. Physical constraints ensure the balance of mass and power, while operational constraints reflect system-specific attributes, such as non-linear power productions of {\ff hydro generators} and water transactions between reservoirs. Figure \ref{fig:Framework} outlines the POMS framework.

Before delving into the model, we first define the general setting and notations of the problem. Our focus is on a {\ff fleet of hydro generators}, owned by a single generation company, over a planning horizon of length $T$. We denote $\mathcal{T}$ and $\mathcal{H}$ as the set of maintenance and operational periods within the planning horizon, respectively. Each maintenance period $t\in\mathcal{T}$ is further subdivided into $H$ operational periods for finer granularity in scheduling.
Each {\ff hydropower system} consists of a reservoir and several {\ff hydro generators} for power production. The sets $\mathcal{R}$ and $\mathcal{G}$ denote the set of reservoirs and {\ff hydro generators}, respectively, while the set $\mathcal{U}$ represents neighboring upstream reservoirs. When necessary, subscript $r$ denotes specific sets associated with the reservoir $r\in\mathcal{R}$. We employ the superscripts $o$ and $f$ to differentiate between operational and maintenance-bound {\ff hydro generators} and those that have failed at the time of planning. Operational {\ff hydro generators} can undergo preventive maintenance, while failed ones are subject to corrective maintenance. The sets $\mathcal{G}^o$ and $\mathcal{G}^f$ denote operational or maintenance-bound and failed {\ff hydro generators}, respectively, at the time of planning. 
\subsection{Objective Function}
The objective is to maximize the expected total profit from {\ff hydropower systems} throughout the planning horizon:
\begin{equation}\label{eq:OBJ} 
\begin{aligned}\max_{\substack{\boldsymbol{e}^+,\boldsymbol{e}^-,\boldsymbol{\alpha}\\\boldsymbol{\beta},\boldsymbol{v},\boldsymbol{z}}}\hspace{2mm}&\sum_{t\in\mathcal{T}}\sum_{h\in\mathcal{H}}\sum_{\omega\in\Omega}p^{\omega}\left( \Gamma_{t,h}^{+,\omega} e_{t,h}^{+,\omega}-\Gamma_{t,h}^{-,\omega} e_{t,h}^{-,\omega}\right)-\sum_{i\in\mathcal{G}}\sum_{t\in\mathcal{T}}\sum_{h \in\mathcal{H}} \sum_{\omega\in\Omega}p^{\omega}\left(C^{su}_{i} \alpha_{i,t,h}^\omega+C^{sd}_{i} \beta_{i,t,h}^\omega\right)\\&-\sum_{i\in\mathcal{G}^o}\sum_{t\in \mathcal{T}} C^{d,t_i^o}_{i,t}v_{i,t}-\sum_{i\in\mathcal{G}^f}\sum_{t \in\mathcal{T}} C^{f}_{i}z_{i,t}\end{aligned}
\end{equation}
\indent
The objective function has four components including the expected revenue from power transactions in the day-ahead market and the costs related to the startup, shutdown, and maintenance of {\ff hydro generators}, to collectively optimize the O\&M planning of {\ff hydropower systems}. The first component captures the expected revenue from power transactions in the day-ahead market, taking into account the stochastic nature of prices for selling excess power ($e^{+,\omega}_{t,h}$), and prices for buying power in case of power scarcity ($e^{-,\omega}_{t,h}$). The terms $\Gamma^{+,\omega}_{t,h}$ and $\Gamma^{-,\omega}_{t,h}$ denote the electricity prices associated with power selling and purchasing, respectively. We incorporate the startup and shutdown costs of {\ff hydro generators} to discourage frequent startups and shutdowns that can lead to increased wear and tear\citep{gagnon2010impact}. We introduce binary variables $\alpha^{\omega}_{i,t,h}$ and $\beta^{\omega}_{i,t,h}$, which indicate if a {\ff hydro generator} $i$ starts up or shuts down during maintenance period $t$ and operational period $h$ in scenario $\omega$, with corresponding startup cost $C^{su}_i$ and shutdown cost $C^{sd}_i$, respectively. Finally, the third and fourth components account for the maintenance costs of operational and failed {\ff hydro generators}. For operational {assets}, the binary variable $v_{i,t}$ denotes the start time of preventive maintenance action; that is, maintenance begins at time $t$ if $v_{i,t}=1$. The dynamic maintenance cost function $C^{d,t^o}_{i,t}$, discussed in Section \ref{sec:DMC}, corresponds to the maintenance decisions of these {assets}. Including this cost function ensures that O\&M scheduling adapts to the latest predictions on failure risks and the remaining life of {\ff hydro generators}. For failed {\ff hydro generators}, $i\in\mathcal{G}^f$, the variable $z_{i,t}$ indicates the start time of the corrective maintenance action. Conducting corrective maintenance actions incurs the corrective maintenance cost $C^f_i$.

\subsection{Constraints}
The objective function (\ref{eq:OBJ}) is subject to the following constraints:
\subsubsection{\textbf{Maintenance Coordination Constraints.}} To ensure high reliability for {\ff hydro generators}, we introduce a dynamic maintenance time limit $\varphi_i^d$ as the first time at which the sensor-driven remaining life probability of the {\ff asset}, denoted as $R_i^{\mathrm{t^o_i}}$, falls below the control threshold $\varrho_i$:  $\varphi_i^d:=min\{t \in \mathcal{T} : P(R_i^{\mathrm{t^o_i}} > t) <\varrho_i\}$. 
The following constraints guarantee that operational {\ff hydro generators} undergo preventive maintenance within this sensor-driven maintenance time window.
\begin{align}\label{eq:m1}
&\sum_{t=1}^{\varphi_i^d}v_{i,t}=1,\quad \forall i \in \mathcal{G}^{\mathrm{o}}
\end{align} 
For failed {\ff hydro generators}, constraints (\ref{eq:m2}) limit the number of corrective maintenance actions to at most one per failed {\ff asset} during the planning horizon $T$. 
\begin{align}\label{eq:m2}
\sum_{t=1}^{T}z_{i,t}\leq1,\quad \forall i \in \mathcal{G}^{\mathrm{f}}
\vspace{-2mm}
\end{align} 
\subsubsection{\textbf{Operations Coordination Constraints.}} 
We categorize the operational constraints of {\ff hydropower systems} into two distinct sets: the first focuses on modeling the operations of {\ff hydro generators} and the second represents the detailed operations of reservoirs. 

\paragraph{\textit{{\ff Hydro generator} Operations Coordination.}}
We represent the startup, shutdown, and commitment decisions of {\ff hydro generators} through binary variables $\boldsymbol{\alpha}$, $\boldsymbol{\beta}$, and $\boldsymbol{x}$, respectively. The water discharge and power output of the {\ff hydro generators} are denoted by continuous variables $\boldsymbol{q}$ and $\boldsymbol{y}$.
    \begin{figure}
    \centering
    \includegraphics[scale=0.5]{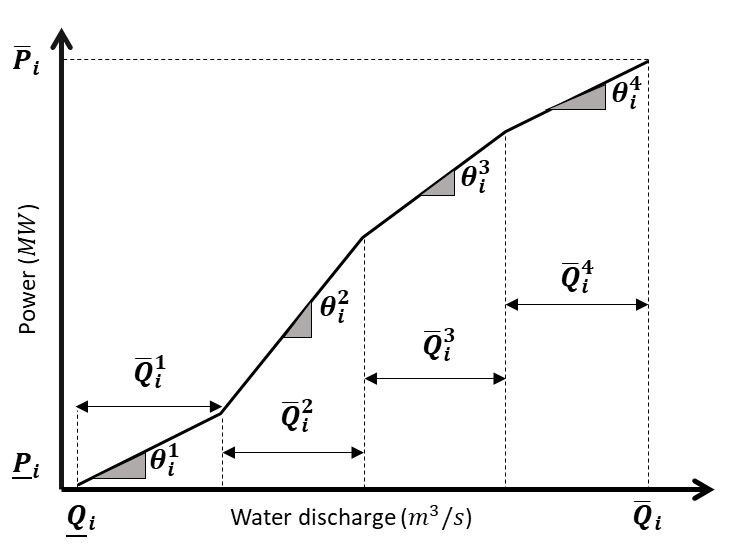}
    \caption{Power production curve of hydro generator $i$}
    \label{fig:powercurve}
\end{figure}
%The power production of a {\ff hydro generator} is a complex function that involves both the kinetic and potential energy of the driving water. These relationships are often represented through non-linear and non-concave power curves\citep{garcia2007risk}. Specifically, 
We describe the power production of each {\ff hydro generator} using an asset-specific power curve to ensure essential fidelity for precise O\&M scheduling of {\ff hydropower systems}. The power output of a {\ff hydro generator} is a non-linear con-concave function of the water discharge rate and the level of stored water in the reservoir. For simplicity, we exclude the head sensitivity from the asset power curves. Additional details on how to include the head sensitivity in power productions can be found in  \citep{conejo2002self,borghetti2008milp}. Constraints \eqref{eq:o1}-\eqref{eq:o5} capture the relationship between the power output of {\ff hydro generators} and the water discharge rate of the assets. 

We approximate the power curves of each {\ff hydro generator} with a piecewise linear function and model the non-linear and non-concave properties of the curve through the use of binary variables. Figure \ref{fig:powercurve} represents a piecewise linear approximation of a generic power curve, estimated with four segments. To formulate the power curve of each {\ff hydro generator $i$}, we define the set $\mathcal{B}_i = {1,...,B_i}$, where $B_i$ represents the number of blocks considered along the power curve of asset $i$. For each block $b\in\mathcal{B}_i$, we define two sets of parameters, $\bar{Q}_i^b$ and $\Theta_i^b$. The parameter $\bar{Q}^b_i$ represents the maximum water discharge capacity (in $m^3/s$) of {\ff hydro generator} $i$ within block $b$, and the parameter $\Theta^b_i$ is the slope of block $b$ of its piecewise linear power curve (in $MW/m^3/s$). For each block $b$, the binary variable $\lambda_{i,t,h}^{b,\omega}$ indicates whether or not the maximum water discharge capacity in block $b$ is reached: $\lambda_{i,t,h}^{b,\omega}=1$ if the water discharged by {\ff hydro generator} $i$ in maintenance period $t$ and operational period $h$ in scenario $\omega$ reaches the discharge limit $\bar{Q}^b_i$ in block $b$. To enforce the water discharge limit $\bar{Q}^b_i$ of the {\ff hydro generators} in each block, we use Constraints \eqref{eq:o1} - \eqref{eq:o3}. 
\begin{align}
&q_{i,t,h}^{1,\omega}\leq \bar{Q}_i^1x_{i,t,h}^{\omega},\hspace{0.5cm}   \forall i \in \mathcal{G}^o,t \in \mathcal{T}, h \in \mathcal{H},\omega\in\Omega\label{eq:o1}\\
&q_{i,t,h}^{b,\omega}\geq \bar{Q}_i^b\lambda_{i,t,h}^{b,\omega},\hspace{0.5cm}\forall i \in \mathcal{G}^{\mathrm{o}}, b \in \mathcal{B}_i, t \in \mathcal{T}, h \in \mathcal{H},\omega\in\Omega\label{eq:o2}\\
&q_{i,t,h}^{b,\omega}\leq \bar{Q}_i^b\lambda_{i,t,h}^{b-1,\omega},\hspace{0.5cm}  \forall i \in \mathcal{G}^o,b \in \mathcal{B}_i, t \in \mathcal{T}, h \in \mathcal{H},\omega\in\Omega\label{eq:o3}
\end{align}
Constraints \eqref{eq:o4} ensure that \begin{inparaenum}[\itshape(1)]
\item if an operational {\ff hydro generator} is committed ($x_{i,t,h}^\omega=1$), its discharged water is at least equal to its minimum water discharge limit $\underline{Q}_i$, and
 \item the discharged water by a committed {\ff hydro generator} in maintenance period $t$ and operation period $h$ in scenario $\omega$ is the sum of the discharged water in each block $b$, i.e., $q_{i,t,h}^{b,\omega}$, plus the minimum water discharge.\end{inparaenum}
\begin{align}
&q_{i,t,h}^{\omega}=\sum_{b=1}^{B_i} q_{i,t,h}^{b,\omega}+\underline{Q}_ix_{i,t,h}^{\omega},\hspace{0.5cm} \forall i \in \mathcal{G}^o, t \in \mathcal{T}, h \in \mathcal{H}\label{eq:o4}
\end{align}
Constraint set \eqref{eq:o5} models the power production of operational {\ff hydro generators} according to their individual piecewise linear production curves. The constraints also enforce the minimum power production limit $\underline{P}_i$ when the {\ff hydro generator} $i$ is committed.
\begin{align}
&y_{i,t,h}^{\omega}=\sum_{b=1}^{B_i} \Theta^b_iq_{i,t,h}^{b,\omega}+\underline{P}_{i}x_{i,t,h}^{\omega},\hspace{0.5cm}\forall i \in \mathcal{G}^{\mathrm{o}}, t \in \mathcal{T}, h \in \mathcal{H}\label{eq:o5}
\end{align}
The last set of constraints establishes the relationship between the startup and shutdown decisions of {\ff hydro generators} and their commitment decisions. We introduce binary variables, $\alpha_{i,t,h}^{\omega}$ and $\beta_{i,t,h}^{\omega}$, to capture the startup and shutdown decisions of {\ff hydro generator} $i$ during maintenance period $t$, operational period $h$, and under scenario $\omega$. Constraint (\ref{eq:o6}) ensures that if a {\ff hydro generator} transitions from not being committed before to being committed, the corresponding startup variable, $\alpha_{i,t,h}^{\omega}$, is set to 1. Similarly, constraint (\ref{eq:o7}) guarantees that if a {\ff hydro generator} transitions from being committed before to not being committed, the shutdown variable, $\beta_{i,t,h}^{\omega}$, is set to 1, indicating its shutdown.
\begin{align}
&\alpha_{i,t,h}^{\omega}+x_{i,t,h-1}^{\omega}-x_{i,t,h}^{\omega}\geq0,\hspace{.5cm}\forall i \in \mathcal{G}^{\mathrm{o}}, t \in \mathcal{T}, h \in \mathcal{H}, \omega\in\Omega\label{eq:o6}\\
&\beta_{i,t,h}^{\omega}-x_{i,t,h-1}^{\omega}+x_{i,t,h}^{\omega}\geq0,\hspace{.5cm}\forall i \in \mathcal{G}^{\mathrm{o}}, t \in \mathcal{T}, h \in \mathcal{H}, \omega\in\Omega\label{eq:o7}
\end{align} 
\paragraph{\textit{Reservoir Operations Coordination.}}
To ensure efficient scheduling and management of hydropower system operations, it is crucial to consider the spatial and temporal relationships between reservoirs. These relationships directly influence water inflows and outflows, which in turn impact the power production of {\ff hydro generators}. In this context, we define the set $\mathcal{U}_r$ as the set of upstream reservoirs connected to reservoir $r$. For any two connected reservoirs $s$ and $r$ from the set $\mathcal{R}$, we define the parameter $\tau_{s,r}$ as the travel time of water flow from upstream reservoir $s$ to downstream reservoir $r$. Furthermore, we denote the volume of reservoir $r$ during maintenance period $t$, operational period $h$, and scenario $\omega$ using the continuous variable $\vartheta_{r,t,h}^{\omega}$.  For each reservoir $r$ at any maintenance period $t$, operational period $h$, and scenario $\omega$, the water inflows to the reservoir consist of the forecasted natural water inflow $F_{i,t,h}^\omega$
%, from a tributary river or snow melt,
as well as the water discharges from immediate upstream reservoirs, $s\in\mathcal{U}_r$. The outflow includes the total discharged water of operational {\ff hydro generators}, $\sum_{i\in\mathcal{G}^o_r}{q_{i,t,h}^{\omega}}$, plus the controlled water spillage $\psi_{r,t,h}^\omega$ from the bypass gates of the reservoir. 

The mass balance constraints for reservoirs are represented by Equations (\ref{eq:o12}). These constraints ensure the conservation of water by relating the water level at operational period $h$ to its level at the previous operational period $h-1$, accounting for the inflows, outflows, and travel time of water between reservoirs. 
\begin{align}\label{eq:o12}
\vartheta_{r,t,h}^{\omega}=&\vartheta_{r,t,h-1}^{\omega}+F_{r,t,h}^{\omega}-K\bigg(\sum_{i\in\mathcal{G}_r}q_{i,t,h}^{\omega}+\psi_{r,t,h}^{\omega}\bigg)+K\sum_{s\in\mathcal{U}_r}\bigg(\sum_{i\in \mathcal{G}_s} q^{\omega}_{i,t,h-\tau_{s,r}}+\psi_{s,t,h-\tau_{s,r}}^{\omega}\bigg)\nonumber\\&\forall r \in \mathcal{R}, t \in \mathcal{T}, h \in \mathcal{H},\omega \in \Omega
\end{align}
Constraint \eqref{eq:o8} ensures that the volume of reservoir $r$ remains within the minimum $\underline{V}_r$ and maximum $\bar{V}_r$ levels, while constraints \eqref{eq:o9} and \eqref{eq:o10} enforce the initial and target water levels of reservoir $r$ at the beginning and end of each maintenance period $t$, respectively. Further details on the initial and target levels of reservoirs will be provided in Section \ref{sec:experiments}.
\begin{align}
&\underline{V}_{r}\leq \vartheta_{r,t,h}^{\omega} \leq \bar{V}_r, \quad\forall r\in\mathcal{R}, t \in \mathcal{T}, h \in \mathcal{H}, \omega \in \Omega\label{eq:o8}\\
&\vartheta_{r,t,0}^{\omega}=V^{0}_{r,t},\hspace{9mm}\forall r\in \mathcal{R}, t \in \mathcal{T}, \omega \in \Omega\label{eq:o9}\\
& \vartheta_{r,t,H}^{\omega}=V^{H}_{r,t},\hspace{9mm}\forall r \in \mathcal{R}, t \in \mathcal{T}, \omega \in \Omega\label{eq:o10}\end{align}
Excessive water inflows or the frequent unavailability of {\ff hydro generators}, either due to failure or maintenance actions, often lead to water spillage in reservoirs. During these occurrences, the storage level of the reservoir approaches its maximum capacity, resulting in water spillage that holds minimal economic value \citep{diniz2008four}. As power production cannot occur during water spillage, it is desirable to minimize such instances \citep{kong2020overview}. To address this, constraint \eqref{eq:o11} is introduced to impose the upper bound $\bar{\Psi}_r$ on the controlled water spillage from reservoir $r$. 
\begin{align}
&\psi_{r,t,h}^{\omega}\leq \bar{\Psi}_r,\hspace{0.5cm}\forall r \in \mathcal{R}, t \in \mathcal{T}, h \in \mathcal{H},\omega \in \Omega\label{eq:o11}
\end{align}
\subsubsection{\textbf{Maintenance and Operations Coupling.}} To establish the interdependence between maintenance decisions and the operations of hydro generators, we introduce  the binary variable $\zeta_{i,t}$, which takes a value of 1 if {\ff hydro generator} $i$ is available to operate during maintenance period $t$, and 0 otherwise.  

The availability of hydro generators is affected by maintenance decisions.  Constraints \eqref{eq:o13} and \eqref{eq:o14} define the availability of {\ff hydro generators} in operational or failed states, respectively. They ensure that if a {\ff hydro generator} $i$ is undergoing preventive maintenance, it will be unavailable for operation in any of the operational periods $h$ during the preventive maintenance action, for all scenarios. The duration of preventive maintenance for hydro generator $i$ is represented by $Y_i^o$.To enforce this availability constraint, we examine whether any preventive maintenance actions have commenced during the preceding maintenance periods ${t-Y_i^o+1,...,t}$.
\begin{align}
&\zeta_{i,t}\leq 1-\sum_{k=0}^{Y^{o}_i-1}v_{i,t-k},\hspace{0.5cm}\forall i \in \mathcal{G}^{\mathrm{o}}, t \in \mathcal{T}, h \in \mathcal{H}, \omega\in\Omega\label{eq:o13}
\end{align}
Constraints \eqref{eq:o14} enforce that a failed {\ff hydro generator} must undergo corrective maintenance for a duration of $Y^{\mathrm{f}}_i$ before it can be available for operation. 
% $Y^{\mathrm{f}}_i$ is the duration of corrective maintenance for the failed {\ff hydro generator} $i$.
\begin{align}
&\zeta_{i,t}\leq \sum_{k=1}^{t-Y^{f}_i}z_{i,k},\hspace{0.5cm}\forall i \in  \mathcal{G}^{\mathrm{f}},t \in \mathcal{T}, h \in \mathcal{H}, \omega\in\Omega\label{eq:o14}
\end{align}
 {\ff Hydro generators} cannot produce power during periods of unavailability and their power production is limited by their maximum $\bar{P}_i$ and minimum $\underline{P}_i$ production capacity. Constraints \eqref{eq:mc3} control the power production of {\ff hydro generators} according to their availability and asset-specific power production capacities.
\begin{equation}
\underline{P}_i\zeta_{i,t}\leq y_{i,t,h}^{\omega}\leq\bar{P}_i\zeta_{i,t},\hspace{0.5cm}\forall i \in  \mathcal{G}, t \in \mathcal{T}, h \in \mathcal{H}, \omega\in\Omega\label{eq:mc3}
\end{equation}
\subsubsection{\textbf{Energy Balance Constraint.}} 
The energy balance constraint \eqref{eq:7} guarantees that the power production and transactions within the hydropower system meet the demand requirements. The power demand is assumed to be known in advance and the surplus or deficit power can be traded in the day-ahead energy market. The constraint ensures that the total produced power by {\ff hydro generators} and purchased power from the market matches the total demand and sold power at any maintenance period $t$, operational period $h$, and in any scenario $\omega$ and is as follows:
\begin{equation}
\sum_{i\in \mathcal{G}^{\mathrm{o}}}y_{i,t,h}^{\omega}+e_{t,h}^{-,\omega}= D_{t,h}+e_{t,h}^{+,\omega},\hspace{0.5cm}  \forall t \in \mathcal{T}, h \in \mathcal{H}, \omega \in \Omega\label{eq:7}
\end{equation}
In summary, the POMS model can be expressed as follows:
\begin{subequations}  \label{eq:APO}
	\begin{align}
\max_{\substack{\boldsymbol{v},\boldsymbol{z},\boldsymbol{q},\boldsymbol{x},\boldsymbol{\lambda},\boldsymbol{y}\\\boldsymbol{\alpha},\boldsymbol{\beta},\boldsymbol{\vartheta},\boldsymbol{\psi},\boldsymbol{\zeta},\boldsymbol{e}}}\hspace{2.5mm}&\hspace{2.5mm} \eqref{eq:OBJ}\\
	\text{s.t.}\hspace{5mm} &\eqref{eq:m1}-\eqref{eq:mc3}\\&
	\{\boldsymbol{q},\boldsymbol{y},\boldsymbol{\psi},\boldsymbol{\vartheta},\boldsymbol{e}\}\in \vv{\mathbb{R}^+}\\&\{\boldsymbol{v},\boldsymbol{z},\boldsymbol{x},\boldsymbol{\lambda},\boldsymbol{\alpha},\boldsymbol{\beta},\boldsymbol{\zeta}\}\in\vv{\{0,1\}}
	\end{align}
	\end{subequations}
The POMS model is computationally demanding due to its size and complexity, which might prevent its practical implementation for some large-scale real-life hydropower systems. In what follows, we propose a tailored solution algorithm to reduce computation time and improve scalability.

\section{Solution Methodology}\label{sec:Algorithm}
To enhance the computational performance of the large-scale POMS model, we employ decomposition techniques to break the POMS model into smaller, more manageable sub-problems. We utilize a combination of Bender's optimality cuts \citep{rahmaniani2017benders} and integer cut to efficiently achieve the optimal O\&M schedules of the original POMS model. 
% This approach allows us to efficiently handle the complexity of the model and obtain optimal hydropower system O\&M schedules within reasonable timeframes.
We start by presenting a two-stage version of the POMS model in a compact form.
\subsection{Two-stage Formulation} 
To enable decomposition-based solution approaches, we reformulate the POMS model \eqref{eq:APO} into a two-stage cost minimization problem.  In this reformulation, the maintenance decisions for the {\ff hydro generators} are determined in the first stage, and the operational decisions given the maintenance decisions, are settled in the second stage. The two-stage stochastic formulation is as follows:
\begin{equation} \label{eq:1st stage}	\min_{\boldsymbol{\nu},\boldsymbol{\zeta}\in\hat{X}}\Bigg\{ \boldsymbol{c}^\top \boldsymbol{\nu} +\mathscr{L}(\boldsymbol{\zeta})\Bigg\}
\end{equation}
where the vectors $\boldsymbol{\nu}$ and $\boldsymbol{\zeta}$ represent the maintenance and availability decisions of {\ff hydro generators}, respectively. The set $\hat{X}$ denotes the feasibility set of the first-stage problem, given by $\hat{X}=\Big\{\boldsymbol{\nu},\boldsymbol{\zeta}|\boldsymbol{A}\hspace{1mm}\boldsymbol{\nu}+B\hspace{1mm}\boldsymbol{\zeta}\leq \boldsymbol{d}\hspace{1mm} \text{and}\hspace{1mm}	 \boldsymbol{\nu},\boldsymbol{\zeta} \in \vv{\{ 0,1 \}}\Big\}$. The term $\mathscr{L}(\boldsymbol{\zeta})$ represents the recourse operations scheduling problem, given the {\ff hydro generator} availability decisions $\boldsymbol{\zeta}$, and is as follows:
\begin{subequations}  \label{eq:2nd stage-relaxd}
	\begin{align}
	\mathscr{L}(\boldsymbol{\zeta})=\min_{\boldsymbol{x},\boldsymbol{y}}\hspace{1mm}&\hspace{1mm}\sum_{\omega\in\Omega}p^{\omega}\big(\boldsymbol{q}^\top\boldsymbol{x}^{\omega}+\boldsymbol{b}^\top \boldsymbol{y}^{\omega}\big)\\ 
	\text{s.t.}&  \hspace{14mm}
	 \boldsymbol{E}\hspace{1mm}\boldsymbol{ x }^{\omega}+\boldsymbol{G}\hspace{1mm}\boldsymbol{y}^{\omega}\leq \boldsymbol{g}-\boldsymbol{D}^{\omega}\boldsymbol{\zeta}\\ 
	&\boldsymbol{x}^{\omega} \in \vv{\{ 0,1 \}}, \boldsymbol{y}^{\omega} \in \vv{\mathbb{R}^+}, \forall \omega \in \Omega\nonumber
	\end{align}
\end{subequations}
Here, the decision vectors $\boldsymbol{x}^\omega$ and $\boldsymbol{y}^\omega$ represent the binary and continuous operational decisions in scenario $\omega\in\Omega$, respectively. 
The presented reformulation allows us to decompose the {\ff hydropower operations scheduling} problem based on maintenance time periods $t\in\mathcal{T}$. More specifically, once the availability decisions of {\ff hydro generators} are determined, the operational decisions of the {\ff hydropower system} at each maintenance period become independent of each other. 
% We represent the decomposed recourse operations scheduling problem by individual maintenance time periods by $\mathscr{L}_t(\boldsymbol{\zeta})$.
% $\mathscr{L}(\boldsymbol{\zeta})=\sum_{t\in\mathcal{T}}\mathscr{L}_t(\boldsymbol{\zeta})$.
% This decomposition strategy has proven to be effective in reducing the computational burden \citep{fallahi2021predictive,okumusoglu2022joint}. 
The \textit{time-decomposed} recourse operations scheduling problem  $\mathscr{L}_t(\boldsymbol{\zeta})$ for each time period $t\in\mathcal{T}$ are as follows:
\begin{subequations}  \label{eq:2nd stage-int}
	\begin{align}
\mathscr{L}_t(\boldsymbol{\zeta})=\min_{\boldsymbol{x}_t,\boldsymbol{y}_t}\hspace{1mm}&\hspace{1mm}\sum_{\omega\in\Omega}p^{\omega}\big(\boldsymbol{q}^\top\boldsymbol{x}_t^{\omega}+\boldsymbol{b}^\top \boldsymbol{y}_t^{\omega}\big)\\ 
	\text{s.t.}&  \hspace{14mm} 
	 \boldsymbol{E}\hspace{1mm}\boldsymbol{ x }^{\omega}_t+\boldsymbol{G}\hspace{1mm}\boldsymbol{y}^{\omega}_t\leq \boldsymbol{g}-\boldsymbol{D}^{\omega}\boldsymbol{\zeta}_t\\ 
    &\boldsymbol{x}_t^{\omega} \in \vv{\{ 0,1 \}}, \boldsymbol{y}_t^{\omega} \in \vv{\mathbb{R}^+}, \forall \omega \in \Omega\nonumber
	\end{align}
\end{subequations}
% The recourse operations scheduling problems, however,  poses integrality restrictions, making it challenging to develop simple solution procedures. 
To develop Bender's algorithm, we relax the binary variables in the second stage to take continuous values between 0 and 1, i.e., $\boldsymbol{x}_t^\omega\in\vv{[0,1]}$ for all $t\in\mathcal{T}$ and $\omega\in\Omega$. We denote the relaxed time-decomposed problem by R-POMS and express it as follows:
\begin{equation} \label{eq:1st stage-R}	\min_{\boldsymbol{\nu},\boldsymbol{\zeta}\in\hat{X}}\Bigg\{ \boldsymbol{c}^\top \boldsymbol{\nu} +\sum_{t\in\mathcal{T}}\bar{\mathscr{L}}_t(\boldsymbol{\zeta})\Bigg\}
\end{equation}
where 
$\bar{\mathscr{L}}_t(\boldsymbol{\zeta})$ represents the relaxed recourse operations scheduling problem.
Solving the relaxed two-stage stochastic model \eqref{eq:1st stage-R} is computationally less demanding than the original model \eqref{eq:1st stage}. However, the R-POMS problem \eqref{eq:1st stage-R} provides a lower bound on the solution of the original POMS model \eqref{eq:1st stage}. Thus, there is a critical need to link the stochastic LP-relaxation model to the original model and recover the true cost as well as O\&M schedules of the original POMS model. In the following sections, we present a two-level decomposition-based algorithm to address the computational and true cost recovery needs.

\subsection{Cost Recovery Algorithm}\label{sec:cost_recovery}
% In this section, we outline the algorithm for retrieving the true cost of the POMS model, i.e., model \eqref{eq:1st stage}, when solving the R-POMS model \eqref{eq:1st stage-R}.
 To establish the link between the POMS and R-POMS model and recover the true schedules and operational costs, we define a set $\bar{\mathcal{A}}_t$ for each maintenance period $t\in\mathcal{T}$. This set contains all possible availability instances of {\ff hydro generators}, resulting in $|\bar{\mathcal{A}}_t|=2^{|\mathcal{G}|}$. We denote an availability instance of {\ff hydro generators} at period $t$ by $\hat{\boldsymbol{\zeta}}_{t}^j\in\bar{\mathcal{A}}_t$. We denote the availability of {\ff hydro generators} in each availability instance by parameters $o_{i,t}^j$ and $n_{i,t}^j$. These parameters are set to 1 if {\ff hydro generator $i$} is not available in the instance $j\in\bar{\mathcal{A}}_t$. If the hydro generator is available, i.e., $\hat{\zeta}_{i,t}^j=1$, the parameters are set to $o_{i,t}^j=-1$ and $n_{i,t}^j=0$. For each maintenance period $t\in\mathcal{T}$ and availability instance $j\in\bar{\mathcal{A}}_t$, we define the binary variable $\mu^j_t$. We enforce the binary variable $\mu^j_t$ to take the value of 1 if availability decision of {\ff hydro generators} at period $t$, i.e., $\boldsymbol{\zeta}_{i,t}$ for all $i\in\mathcal{G}$, corresponds to the $j^{th}$ availability instance in set $\bar{\mathcal{A}}_t$, i.e., $\hat{\zeta}^j_t$.  The following constraint enforces the relationship between the availability instance and availability decisions of {\ff hydrogenerators}:

 %Recall that the binary decision vector $\boldsymbol{\zeta}_{t}$ represents availability decisions for {\ff hydro generators} at period $t$, where its $i^\text{th}$ element is equal to $1$ if the {\ff hydro generator} $i$ is available to operate at that time and $0$ otherwise. Once the availability decisions for {\ff hydro generators} are determined, the true cost of the recourse functions, $\sum_{t\in\mathcal{T}}\mathscr{L}(\boldsymbol{\zeta})$, becomes computable. In other words, the availability decisions provide sufficient information to calculate the cost of the R-POMS \eqref{eq:2nd stage-int} model for every $t$. 

% To establish the link and recover the true cost, for each period $t\in\mathcal{T}$, we define a set $\mathcal{A}_{t}$ containing all possible availability instances of {\ff hydro generators}, i.e.,  $|\mathcal{A}_t|=2^{|\mathcal{G}|}$. For each period $t\in\mathcal{T}$, we denote an availability instance of {\ff hydro generators} by $\hat{\boldsymbol{\zeta}}_{t}^j\in\mathcal{A}_t$. We introduce the binary variable $\mu_t^j$, which takes the value 1 if the availability solution $\boldsymbol{\zeta}_t$ corresponds to the specific availability status $\hat{\boldsymbol{\zeta}}_{t}^j$, and 0 otherwise. We enforce the following constraints, referred to as \textit{cost recovery cuts}, to establish such relationships:

\begin{equation}\label{eq:integercuts}
    \mu^j_t\geq\sum_{i\in \mathcal{G}}\bigg(o_{i,t}^j\cdot(1-\zeta_{i,t})-n_{i,t}^j\bigg)+1,\quad \forall t \in \mathcal{T}
\end{equation}
%The binary variable $\mu_t^j$ takes the value 1 if the value of availability decisions $\boldsymbol{\zeta}_{i,t}$ for all $i\in\mathcal{G}$ corresponds to the availability status $\hat{\boldsymbol{\zeta}}_{t}^j$, and 0 otherwise.
% This constraint forces the binary variable $\mu^j_t$ take the value of 1 if {\ff hydro generators} availability decisions at period $t$, i.e., $\boldsymbol{\zeta}_{i,t}$ for all $i\in\mathcal{G}$, corresponds to the $j^{th}$ availability instance in set $\mathcal{A}_t$, i.e., $\hat{\zeta}^j_t$.To enable this relationship, for each period $t\in\mathcal{T}$, 
% The parameters $o_{i,t}^j$ and $n_{i,t}^j$ are set to 1 when the status of {\ff hydro generator $i$} is non-available in the availability instance $j\in\mathcal{A}_t$. If the hydro generator is available, i.e., $\hat{\zeta}_{i,t}^j=1$, the parameters are set to $o_{i,t}^j=-1$ and $n_{i,t}^j=0$. 
 We refer to this constraint set as \textit{integer cuts} and express it in its compact form as $\mu_t^j\geq\boldsymbol{u}_t^j-(\boldsymbol{q}_t^j)^\top\boldsymbol{\zeta}_t$. We define the parameter  $\delta^j_t$ as the cost difference between the original and relaxed operational recourse problem. More precisely, $\delta^j_t= \mathscr{L}_t(\hat{\boldsymbol{\zeta}}_{t}^j)-\bar{\mathscr{L}}_t(\hat{\boldsymbol{\zeta}}_{t}^j)$, for any solution $\boldsymbol{\zeta}_{t}$ that implies $\hat{\boldsymbol{\zeta}}_{t}^j$. Consequently, the following equation holds for any availability decisions $\boldsymbol{\zeta}_t$ and maintenance period $t\in\mathcal{T}$:
\begin{equation}
\mathscr{L}_t(\boldsymbol{\zeta}_t)=\bar{\mathscr{L}}_t(\boldsymbol{\zeta}_t)+\min_{\boldsymbol{\mu}}\Big\{\sum_{j\in\mathcal{A}_t}\delta_t^j\mu_t^j:\mu_t^j\geq\boldsymbol{u}_t^j-(\boldsymbol{q}_t^j)^\top\boldsymbol{\zeta}_t\Big\}
\end{equation}
% The optimal solution of this problem is clear given the fact that for only one $j\in\mathcal{A}_j$ is the binary variable $\mu_t^j$ forced to be 1, i.e., $\boldsymbol{u}_t^j-(\boldsymbol{q}_t^j)^\top\boldsymbol{\zeta}_t=1$.
We represent the optimal solution of this problem by $\mathscr{L}_t(\boldsymbol{\zeta}_t)=\bar{\mathscr{L}}_t(\boldsymbol{\zeta}_t)+\delta^{j^*}_t$ in which $j^*$ represents the index of availability statuses implied from the availability decisions $\boldsymbol{\zeta}_t$. By introducing these integer cuts to the R-POMS model, the following model can attain the optimal objective cost and decisions of the original POMS model while solving the relaxed recourse function $\bar{\mathscr{L}}_t(\boldsymbol{\zeta})$.
\begin{subequations}  \label{eq:Final model}
	\begin{align}
 \hspace{4pt}	\min_{\boldsymbol{\nu},\boldsymbol{\zeta},\boldsymbol{\mu}} \hspace{4mm}&\boldsymbol{a}^\top \boldsymbol{\nu}+\sum_{t\in\mathcal{T}}\boldsymbol{\delta}_t^\top\boldsymbol{\mu}_t +\sum_{t\in\mathcal{T}}\bar{\mathscr{L}}_t(\boldsymbol{\zeta})\\	\text{s.t.}  &\hspace{1.5mm}\boldsymbol{A}\hspace{1mm}\boldsymbol{\nu}+B\hspace{1mm}\boldsymbol{\zeta}\hspace{18 mm}\leq \boldsymbol{d}\\&
 \mu_t^j\geq\boldsymbol{u}_t^j-(\boldsymbol{q}_t^j)^\top\boldsymbol{\zeta}, \quad \forall t \in \mathcal{T}, \forall j \in \bar{\mathcal{A}}_t\\&
	\quad{ \boldsymbol{\nu},\boldsymbol{\zeta} \in \vv{\{ 0,1 \}}},\boldsymbol{\mu}\in\{0,1\}^{T\times2^{|\mathcal{G}|}}\nonumber
	\end{align}
\end{subequations}
 % which can attain the optimal objective cost and decisions of the original model while encompassing the relaxed recourse function $\bar{\mathscr{L}}_t(\boldsymbol{\zeta})$. 
Note that, in most cases, incorporating a subset of the complete availability set $\bar{\mathcal{A}}_t$ appears to be sufficient for recovering the true cost since the maintenance cost minimizer $\boldsymbol{a}^\top \boldsymbol{\nu}$ plays the dominant role in determining the optimal maintenance time windows. 
% We initialize the availability sets of {\ff hydro generators} $\mathcal{A}_t$ with some potential availability instances $\hat{\boldsymbol{\zeta}}^j_t$ and record their corresponding true operational costs, i.e., $\mathscr{L}(\hat{\boldsymbol{\zeta}}^j_t)$. 
In addition, the cost difference between the relaxed R-POMS model and the original POMS model is typically small. Therefore, in practice, to recover the true cost and solutions of the POMS model, it is sufficient to check a few other maintenance time windows. In the next section, we incorporate the developed cost recovery algorithm into a two-level solution algorithm to speed up the solution process of the POMS model.

\subsection{Two-level Solution Algorithm}
We use a combination of Benders optimality cuts and integer cuts to develop a two-level solution algorithm. The main idea is to first solve the relaxed R-POMS by utilizing Benders optimality cuts and then, add integer cuts, if needed, to retrieve the true operational cost of {\ff hydropower systems} and repeat this process until optimal convergence. The Benders cuts approximate the operational costs of the R-POMS model $\bar{\mathscr{L}}_t(\boldsymbol{\zeta})$ through supporting hyperplanes, and cost recovery cuts add the cost difference between the relaxed and the original model, i.e., $\sum_{t\in\mathcal{T}}(\mathscr{L}_t(\boldsymbol{\zeta})-\bar{\mathscr{L}}_t(\boldsymbol{\zeta}))$, back to the total expected O\&M cost the R-POMS model through integer cuts. The pseudo-code of the developed two-level solution algorithm is provided in the appendix.

We start elaborating on the algorithm by defining some key notations and definitions. We denote the tolerance level and convergence flag of the Benders and cost recovery algorithms by $\epsilon^\mathcal{B}\geq0$ and $\mathcal{F}^\mathcal{B}$ and $\epsilon^\mathcal{C}\geq0$ and $\mathcal{F}^\mathcal{C}$, respectively. We set the convergence flags to \textit{False} at the very beginning. We define the iteration number of the Benders and cost recovery algorithms by $\ell^\mathcal{B}$ and $\ell^\mathcal{C}$, respectively. The set $\mathcal{BD}$ represents the set of generated Bender's optimality cuts with $k\in\mathcal{BD}$ representing the cuts index number. We set up the master problem as follows: 
 \begin{subequations}  \label{model:master}
	\begin{align}
 \hspace{4pt}	\min_{\boldsymbol{\nu},\boldsymbol{\zeta},\boldsymbol{\mu},\boldsymbol{\chi}} \hspace{4mm}&\boldsymbol{c}^\top \boldsymbol{\nu}+\sum_{t\in\mathcal{T}}\boldsymbol{\delta}_t^\top\boldsymbol{\mu}_t + \sum_{t\in\mathcal{T}}\boldsymbol{\chi}_t\\	\text{s.t.}  &\hspace{1.5mm}\boldsymbol{A}\hspace{1mm}\boldsymbol{\nu}+B\hspace{1mm}\boldsymbol{\zeta}_t\leq \boldsymbol{d}\label{eq:benders-cons1}\\& \boldsymbol{\chi}_t\geq\boldsymbol{\alpha}_{t}^{k}-\boldsymbol{\beta}^{k}_t\boldsymbol{\zeta}_t, \quad \forall t \in \mathcal{T},k\in\mathcal{BD}\label{eq:optimcut}\\&
 \mu_t^j\geq\boldsymbol{u}_t^j-(\boldsymbol{q}_t^j)^\top\boldsymbol{\zeta}_t, \quad \forall t \in \mathcal{T}, \forall j \in \mathcal{A}_t\label{eq:integercuts-compact}\\&
	\quad{ \boldsymbol{\nu},\boldsymbol{\zeta}\in \vv{\{ 0,1 \}}},\boldsymbol{\mu}\in\{0,1\}^{|\mathcal{A}|},\boldsymbol{\chi}\in \vv{\mathbb{R}}\nonumber
	\end{align}
\end{subequations}
in which the free variable $\boldsymbol{\chi}_t$ represents the cost of the relaxed operational problem $\bar{\mathscr{L}}_t(\boldsymbol{\zeta}_t)$ for each $t\in\mathcal{T}$. Constraints \eqref{eq:optimcut} and \eqref{eq:integercuts-compact} are the Benders optimality cuts and integer cuts, respectively.
Let $(\boldsymbol{\nu}^{\ell^\mathcal{B}},\boldsymbol{\zeta}^{\ell^\mathcal{B}},\boldsymbol{\chi}^{\ell^\mathcal{B}})$ be the optimal solution of this problem at iteration $\ell^\mathcal{B}$. If no constraint  \eqref{eq:optimcut} is present, $\boldsymbol{\chi}^{\ell^\mathcal{B}}$ is set equal to $-\infty$ and is not considered in the computation of $(\boldsymbol{\nu}^{\ell^\mathcal{B}},\boldsymbol{\zeta}^{\ell^\mathcal{B}})$.  Given the obtained {\ff hydro generators} availability schedules $\boldsymbol{\zeta}_t^{\ell^\mathcal{B}}$, for each $t\in\mathcal{T}$ we solve the relaxed operational sub-problems $\bar{\mathscr{L}}_t(.)$:
\begin{subequations}  \label{eq:benders-sub}
	\begin{align}
\bar{\mathscr{L}}_t(\boldsymbol{\zeta}_t^{\ell^\mathcal{B}})=\min_{\boldsymbol{x}_t,\boldsymbol{y}_t}\hspace{1mm}&\hspace{1mm}\sum_{\omega\in\Omega}p^{\omega}\big(\boldsymbol{q}^\top\boldsymbol{x}_t^{\omega}+\boldsymbol{b}^\top \boldsymbol{y}_t^{\omega}\big)\\ 
	\text{s.t.}&  \hspace{10mm} \label{eq:2nd stage-relaxd-cons}
	 \boldsymbol{E}\hspace{1mm}\boldsymbol{ x }^{\omega}_t+\boldsymbol{G\hspace{1mm}y}^{\omega}_t\leq \boldsymbol{g}-\boldsymbol{D}^{\omega}\boldsymbol{\zeta}^{\ell^\mathcal{B}}_t\\ 
    &\hspace{5mm}\boldsymbol{x}_t^{\omega} \in \vv{[ 0,1 ]}, \boldsymbol{y}_t^{\omega} \in \vv{\mathbb{R}^+}, \forall \omega \in \Omega\nonumber
	\end{align}
\end{subequations}
We represent the dual multipliers associated with the optimal solution of the master problem at iteration $\ell^\mathcal{B}$
 by $\boldsymbol{\pi}^{\ell^\mathcal{B}}_{t}$ for each $t\in\mathcal{T}$. If the inequality
 \begin{equation}\boldsymbol{\chi}_t^{\ell^\mathcal{B}}<\sum_{\omega\in\Omega}p^{\omega}{(\boldsymbol{\pi}_{t}^{\omega,\ell^\mathcal{B}})}^\top\boldsymbol{g}-\big(\sum_{\omega\in\Omega}p^{\omega} {(\boldsymbol{\pi}_{t}^{\omega,\ell^\mathcal{B}})}^\top\boldsymbol{D}^\omega\big)\boldsymbol{\zeta}_t^{\ell^\mathcal{B}}\end{equation}
 we set $k=k+1$ and define the Benders multipliers as: 
\begin{equation}
\boldsymbol{\alpha}_{t}^{k+1}=\sum_{\omega\in\Omega} {(\boldsymbol{\pi}_{t}^{\omega,\ell^\mathcal{B}})}^\top\boldsymbol{g}\hspace{0.5cm},\hspace{0.5cm}
\boldsymbol{\beta}^{k
+1}_t=\sum_{\omega\in\Omega} {(\boldsymbol{\pi}_{t}^{\omega,\ell^\mathcal{B}})}^\top\boldsymbol{D}^\omega
\end{equation}
We generate and add the new Benders optimality cuts to the master problem and solve it again. We repeat this process until the inequality $\sum_{t\in\mathcal{T}}\bar{\mathscr{L}}_t(\boldsymbol{\zeta}_t^{\ell^\mathcal{B}})\leq(1+\epsilon^\mathcal{B})\sum_{t\in\mathcal{T}}\boldsymbol{\chi}_t^{\ell^\mathcal{B}}$proves to be true, i.e., the Benders convergence criteria are satisfied.
We set the Benders convergence flag to $\mathcal{F}^\mathcal{B}=True$ and denote the optimal
% solution of the R-POMS model by $({\boldsymbol{\nu}^\mathcal{B}}^{*},{\boldsymbol{\zeta}^\mathcal{B}}^{*},{\boldsymbol{\chi}^\mathcal{B}}^{*})$,and optimal 
cost by $\rho^{\ell^\mathcal{B}}$. 
Next, we check to see if the true operational cost of the original POMS model differs significantly from the operational cost $\rho^{\ell^\mathcal{B}}$ obtained via the R-POMS model and take appropriate actions. 

For each $t\in\mathcal{T}$, we first check whether the availability set $\mathcal{A}_t$ contains the availability solution $\boldsymbol{\zeta}_t^{\ell^\mathcal{B}}$: If it does not, we update the set of availability instances accordingly, i.e., $\mathcal{A}_t=\mathcal{A}_t\cup\{\boldsymbol{\zeta}_t^{\ell^\mathcal{B}}\}$ for every $t\in\mathcal{T}$. We set $j=j+1$, and solve the sub-problems \eqref{eq:2nd stage-int} to evaluate the operational cost difference $\delta^j_t$ when the availability status of {\ff hydro generators} is equivalent to the $\boldsymbol{\zeta}_t^{\ell^\mathcal{B}}$. The total cost difference is equal to $\sum_{t\in\mathcal{T}}\delta^j_t=\sum_{t\in\mathcal{T}}\big(\mathscr{L}_t(\boldsymbol{\zeta}_t^{\ell^\mathcal{B}})-\bar{\mathscr{L}}_t(\boldsymbol{\zeta}_t^{\ell^\mathcal{B}})\big)$. Given the total operational cost difference between the R-POMS and POMS model
, we evaluate the $\epsilon^{\mathcal{C}}$ optimality criterion of the cost recovery algorithm as follows:
\begin{equation}\label{eq:CostRecovery_Optimality}
     \rho^{\ell^\mathcal{B}}+\sum_{t\in\mathcal{T}}\delta^j_t\leq(1+\epsilon^\mathcal{C}) \rho^{\ell^\mathcal{B}}
\end{equation}
If this inequality holds true, the optimal solution is achieved, and we set the optimality flag $\mathcal{F}^\mathcal{C}=True$. In other words, the optimal cost of the model in Equation  \eqref{model:master} proves to be close enough to its corresponding true cost, i.e., the optimal cost of the POMS model \eqref{eq:1st stage}. Otherwise, we add new binary variables $\mu_t^j$ along with their corresponding costs $\delta_t^j$, and the newly generated cost recovery cuts $ \mu_t^j\geq\boldsymbol{u}_t^j-(\boldsymbol{q}_t^j)^\top\boldsymbol{\zeta}_t$ to the model in  \eqref{model:master}. We set the convergence flag of the Benders algorithm $\mathcal{F}^\mathcal{B}$ to \textit{False}, and solve the master problem in \eqref{model:master} again. We repeat this process until both Benders and cost recovery convergence criteria are satisfied.

Note that the two-stage model has complete recourse and there is no need to generate feasibility cuts in the Benders algorithm. Note also that by incorporating a subset of hydrogenerator availability instances and its corresponding costs from the beginning through the constraints \eqref{eq:integercuts-compact} we may not need to run the cost recovery algorithm in practice.

\section{Experiments}\label{sec:experiments}
In this section, we develop a comprehensive experimental platform that uses rolling horizon simulation to demonstrate the effectiveness of the POMS framework in scheduling O\&M decisions for a fleet of hydro generators. {{We leverage real-world hydropower data to assess and validate the practicality of the developed sensor-driven O\&M scheduling framework. We use the vibration signals of bearings in mils within a hydropower system setting. The x-axis indicates operational time in weeks, and the y-axis denotes the vibration displacement in mils. Consistent with International Organization for Standardization (ISO) recommendations, we set the failure threshold of bearings to 10 mil. Consequently, we denote the failure time as the first time that the vibration signals exceed this failure threshold. During preprocessing the vibration signals, we address sensor noise and errors by removing outliers, interpolating missing values to maintain data consistency, and applying a moving average technique to emphasize trends and patterns that signify bearing wear and potential failure. 
% This dataset, which provides a comprehensive overview of bearing degradation via vibration displacement measurements, was collected by HRI to enable digital transformation in the hydropower power industry through projects such as sensor-driven predictive maintenance scheduling of hydropower systems.
 The vibration data has been used to estimate the degradation signal parameters defined in degradation formulation \eqref{eq:Deg}.}} %The parameters are then used to emulate degradation of generators in our framework.

As our benchmark, we use industry-standard time-based maintenance (TBM) policy. Unlike our sensor-driven approach that employs real-time data, the TBM policy schedules maintenance activities periodically, irrespective of sensor inputs. In what follows, we provide details about the test system used in our experiments. We next outline the experimental framework and introduce the performance metrics. Finally, we present the outcomes of the experiments. Detailed explanations of the input datasets are available in the Appendix.
% We test the performance of the POMS and TBM models in a rolling horizon fashion and capture key performance metrics. 
% In the following subsections, we first describe the test system and associated data sets. Then, we outline the experimental framework and introduce the performance metrics. Finally, we present the outcomes of the experiments.
\begin{table*}[!htbp]
\caption{Reservoir Characteristics}
\centering
\adjustbox{width=0.7\textwidth}{% NOTE: N = 10 in this table.
\begin{tabular}{ l C{1cm} C{1cm} C{1cm} C{1.2cm}C{1cm} C{1cm} C{1cm} C{1.2cm}} 
\toprule
    \textbf{Reservoirs ($r$)}& \textbf{1}& \textbf{2}& \textbf{3}& \textbf{4}& \textbf{5}& \textbf{6}& \textbf{7}& \textbf{8}\\ \midrule \midrule
	\textbf{HydroG}	
	&1, 2 &3, 4 &5, 6&7, 8&9, 10&11, 12&13, 14&15, 16\\
% 	\\ \\\\[-3.8\medskipamount] 
	\hline 
% 	\\\\[-3.8\medskipamount]
	$\boldsymbol{\underline{V}_r}$ ($Hm^3$)	&6 &6 &6 &6 &6 &6&6 &6  \\
% 	\\\\[-3.8\medskipamount]
	\hline
	% \\\\[0.0001\medskipamount]
	\multirow{1}{*}{$\boldsymbol{\bar{V}_r}$ ($Hm^3$)}	
	&225 &162 &1200 &66 &26 &2586&115 &181\\\hline \multirow{1}{*}{$\boldsymbol{V}_{r}^0$ ($Hm^3$)}	
	&100 &80 &790 &33 &13 &1200&50&90\\\hline \multirow{1}{*}{$\boldsymbol{\bar{F}}_{r}$ ($Hm^3$)}	
	&0.051 &0.058 &0.603 &0.051 &0.051 &0.199&0.500 &0.048\\ \bottomrule
\end{tabular}}\label{table:Rsrvrs}
\end{table*}

    \begin{figure}[!htbp]
    \centering
    \includegraphics[scale=0.55]{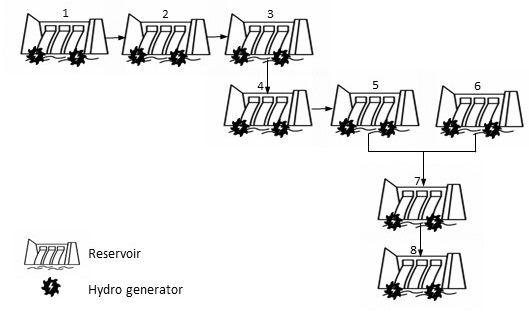}
    \caption{Spatial interconnection of reservoirs} \label{fig:Casestudy}
\end{figure}  
\subsection{Test System and Data}\label{sec:experiments-A}
   We study a 24-week joint O\&M planning of {\ff hydropower systems} with weekly  maintenance and hourly operational decisions, respectively. 
   % Hence, $T=24$ and $H=168$.
   The test system under consideration is modified from a realistic test case introduced in \citep{conejo2002self}. This setup consists of eight reservoirs with different power generation and storage capacities. The schematic representation of the connections between reservoirs is shown in Figure \ref{fig:Casestudy}. Each reservoir is equipped with two {\ff hydro generators}. For any two connected reservoirs denoted as $s$ and $r$ (where $s, r \in \mathcal{R}$), a one-hour time delay (traveling time for water) is assumed. We denote the temporal lag between reservoirs by $\tau_{s,r}=1$. We outline the characteristics of reservoirs and their corresponding {\ff hydro generators} in Tables \ref{table:Rsrvrs} and \ref{table:HGs}, respectively. The non-linear power production curves of {\ff hydro generators} are approximated using four equidistant linear segments. 
The parameters of these approximated power curves, including slopes, minimum and maximum power outputs, and corresponding water discharges, are listed in Table \ref{table:HGs}. We set the preventive (\(C^p\)) and corrective (\(C^f\)) maintenance costs for {\ff hydro generators} at \$200,000 and \$800,000, respectively. Crew deployment costs are omitted (\(C^v=0\)). The duration for preventive maintenance is one week (\(Y^p=1\)), while corrective maintenance takes twice as long (\(Y^f=2Y^p\)) to account for the additional coordination needed for unexpected failures.  %We conduct all experiments using the Gurobi 10 optimizer. 
%**********************************************************
\begin{table*}[!htbp]
\caption{Hydrogenerator Characteristics}
\centering
\adjustbox{width=0.7\textwidth}{% NOTE: N = 10 in this table.
\begin{tabular}{ l C{1.1cm} C{1.1cm} C{1.1 cm} C{1.45 cm} C{1.45 cm} C{0.75 cm} C{0.7 cm} C{0.7 cm}C{0.7 cm}} 
\toprule
    \textbf{HydroG} &$\boldsymbol{\underline{P}_i}$ \scriptsize{(MW)}& $\boldsymbol{\bar{P}_i}$ \scriptsize{(MW)}&$\boldsymbol{C_i^{su}}$ & $\boldsymbol{\underline{Q}_i}$ \scriptsize{(MW/$m^3/s$}) & $\boldsymbol{\bar{Q}_i^{b}}$ \scriptsize{(MW/$m^3/s$})& $\boldsymbol{\Theta^1_i}$& $\boldsymbol{\Theta^2_i}$&$\boldsymbol{\Theta^3_i}$&$\boldsymbol{\Theta^4_i}$ \\ \midrule \midrule
% 	\multirow{1}{*}{\textbf{Reservoir}}
% 	 &1& 2& 3&4 & 5&6&7&8 \\ \\\\[-3.8\medskipamount] \hdashline \\\\[-3.8\medskipamount]
	\multirow{1}{*}{\textbf{1, 2}}	
	& 0.81 & 14.31 & \$55 &1 &7.5  &0.9  &0.4&0.3&0.2 \\
% 	\\\\[-3.8\medskipamount]
	\hline
% 	\\\\[-3.8\medskipamount]
	\multirow{1}{*}{\textbf{3, 4}}	& 1.185 & 34.76 & \$75 & 2.5 &19.75  &0.5&0.4&0.6&0.2 \\
% 	\\\\[-3.8\medskipamount]
	\hline 
% 	\\\\[-3.8\medskipamount]
	\multirow{1}{*}{\textbf{5, 6}}	
	& 2.025 & 69.53 & \$100 &7 & 56.25&0.3 &0.2&0.4&0.3 \\
% 	\\\\[-3.8\medskipamount]
	\hline
% 	\\\\[-3.8\medskipamount]
	\multirow{1}{*}{\textbf{7, 8}}	& 1.929 & 58.19 & \$125 &9.5  & 80.375 & 0.2&0.2&0.15&0.15
	\\
% 	\\\\[-3.8\medskipamount] 
	\hline 
% 	\\\\[-3.8\medskipamount]
	\multirow{1}{*}{\textbf{9, 10}}	
	& 1.83 & 93.33 & \$175 & 9 &76.25 &0.2  &0.5&0.3&0.2 \\
% 	\\\\[-3.8\medskipamount]
	\hline
% 	\\\\[-3.8\medskipamount]
	\multirow{1}{*}{\textbf{11, 12}}	& 9.765 & 416.64 & \$750 & 7 & 58.125&1.4&3.1&1.6&0.9 \\
% 	\\\\[-3.8\medskipamount]
	\hline 
% 	\\\\[-3.8\medskipamount]
	\multirow{1}{*}{\textbf{13, 14}}
	& 12.189 & 579.82 & \$1000 &14.5 &119.5&0.85  & 1.6&1.3&1\\ 
% 	\\\\[-3.8\medskipamount]
	\hline 
% 	\\\\[-3.8\medskipamount]
	\multirow{1}{*}{\textbf{15, 16}}& 13.473 & 275.45 & \$500 &15&124.75  &0.9 &0.4  &0.6&0.2\\  \bottomrule
\end{tabular}}\label{table:HGs}
\end{table*}

\subsection{Experimental Framework and Benchmark Policy}\label{sec:experiments-B}

Experimental framework builds on a rolling horizon approach that follows iterative applications of O\&M optimization and simulation. The framework starts with solving the 24-week stochastic joint O\&M scheduling problem (\textit{planning step}). We then fix the optimal {\ff hydro generator} maintenance schedules for the initial eight weeks, known as the \textit{freeze period}. We then simulate the sequence of operational and maintenance-related events that happen following the scheduled maintenance actions (\textit{execution step}). Specifically, the execution step conducts the following tasks: \begin{inparaenum}[\itshape (1)]\item  For operational {\ff hydro generators} at the planning step, we evaluate if they successfully undergo planned preventive maintenance or encounter unexpected failures. Failures are identified by examining the magnitude of corresponding degradation signals and failure thresholds. If a {\ff hydro generator} fails prior to its scheduled preventive maintenance, it remains offline for the rest of the freeze period. Otherwise, the {\ff hydro generator} undergoes preventive maintenance for one week. \item If any failed {\ff hydro generators} are detected at the planning step, we check whether they are scheduled for corrective maintenance within the freeze period.  If so, the {\ff hydro generator} receives corrective maintenance for two weeks\end{inparaenum}.  

% and becomes available for the rest of the freeze period
After the execution step, we resolve the operations scheduling model 
% for each week individually
using the derived availability of {\ff hydro generators} and revealed market price and water inflow scenarios. 
During the execution step, we track key performance indicators such as the number of preventive and corrective maintenance actions, number of failures, power productions, net power transactions, and expenditures. We also keep track of resource utilization metrics, specifically the availability rate and unused life of {\ff hydro generators} as well as water spillage of reservoirs. 
After the execution step ends, we roll forward 8 weeks to plan the subsequent 24-week O\&M scheduling for the hydropower systems. We assign new degradation signals randomly to repaired {\ff hydro generators} from the database, emulating the degradation of brand-new assets. Dynamic maintenance costs for each {\ff hydro generator} are updated based on the observed degradation sensor data during freeze periods. This process is executed eight times to cover a 64-week period. We repeat this entire process ten times with different degradation signals and initial ages for 16 {\ff hydro generators}. The reported metrics are the average of collected metrics from these ten individual experiments.

\subsection{Experimental Results}
In this section, we present and discuss the experimental results.

\subsubsection{\textbf{Experiment I-Benchmark Analysis:}} \label{sec:Results_Benchmark}
In this experiment, we showcase the efficacy of POMS and TBM models in the O\&M scheduling of {\ff hydro generators}. Table \ref{table:benchmark} presents key performance metrics of both models, including preventive and corrective maintenance actions, failures, outages, and economic indicators such as operational revenue and net profit. 
\begin{table}[!htbp] \centering
	\caption{Benchmark Analysis}
	\adjustbox{width=0.5\textwidth}{%
	\begin{tabular}{l C{2cm} C{2cm}} 
        \toprule \textbf{Performance Metric} & \textbf{ TBM } & \textbf{ POMS }  \\ \midrule \midrule
		%\textbf{Metrics }&& & \textbf{N=15} & \textbf{N=15 } \\
	    \textbf{Preventive Actions} & 15& 12.3\\
	    \textbf{Corrective Actions}  & 6.3 &1.7\\
	     \textbf{Failures}  & 6.7& 1.9\\
	    \textbf{Outages}  & 21.7& 14.2\\ \midrule
	    		\textbf{Availability}  & 95.25\%& 97.83\%\\
 		\textbf{Unused Life (Weeks)}  & 59.3& 21.14\\
 						\textbf{Spilled Water ($Hm^3$)}  & 87.85 & 27.01\\
		\textbf{Power Production ($MW$)} & 12.22 \textbf{M} & 12.30 \textbf{M} \\\textbf{Net Power Transaction ($MW$)} & 8.08 \textbf{M} & 8.16 \textbf{M} \\
 \midrule
		\textbf{Maintenance Cost} & \$8.04 \textbf{M}& \$3.82 \textbf{M}\\
		\textbf{Operational Revenue}& \$190.23 \textbf{M} & \$194.06 \textbf{M}\\
		\textbf{Net Profit} & \$182.19 \textbf{M} & \$190.24 \textbf{M}\\
		\bottomrule
	\end{tabular}}\label{table:benchmark}
\end{table}

 The POMS model capitalizes on sensor data to identify and predict critical conditions in {\ff hydro generators} and schedule necessary maintenance actions. This results in an 18\% reduction in the number of preventive maintenance actions compared to the TBM model, which schedules such actions at periodic intervals without consideration of the health conditions of {\ff hydro generators}. Another important metric is the \textit{unused life} of {\ff hydro generators}. The unused life of an asset refers to the time between the planned maintenance time and the potential failure time under no maintenance treatment. The POMS model is able to decrease the unused life of the {\ff hydro generators} 64.36\%.  By utilizing the sensor-driven dynamic maintenance costs, the POMS model penalizes poorly-timed maintenance actions based on the predicted failure risks of the {\ff hydro generators}. The lower number of failure occurrences in the POMS model demonstrates the optimal timing of preventive maintenance actions which results in a 34.56\% decrease in the total number of {\ff hydro generator} outages. Consequently, the POMS model improves asset availability by 2.71\% compared to the TBM model. 
 
 The enhanced availability of {\ff hydro generators} under the POMS model enables more efficient power generation and reduces water spillage with no economic value by 69.25\%. It also results in a 0.65\% increase in total power production. The higher availability rate in the POMS model directly impacts the power transactions capability in the day-ahead energy market, i.e., the ability of hydropower systems to meet demand and capitalize on high market prices by selling the surplus power. Specifically, the net power transactions in the POMS model are increased by 0.99\% in comparison to the TBM model. The total \textit{net power transactions} in the day ahead energy market is defined as the difference between the sold surplus power and purchased power to substitute production deficiency over time. From a financial standpoint, the POMS model decreases the total maintenance costs by 52.49\%, leading to a 4.41\% improvement in net profit when compared to the TBM model.
\subsubsection{\textbf{Experiment II: Maintenance Duration Effect}}\label{sec:Results_Duration}
This experiment aims to investigate the impact of varying maintenance durations on the performance metrics of POMS and TBM models. Specifically, we examine how different ratios of corrective to preventive maintenance durations \(\left(\frac{Y^f}{Y^p}\right)\) affect the efficacy of O\&M schedules in both models. This variable ratio addresses real-world scenarios where maintenance durations can vary due to factors such as the severity of failures and staff availability. We alter the \(\frac{Y^f}{Y^p}\) ratio, ranging from 1 to 3 while keeping all other model parameters constant. The performance metrics of both models for each ratio are summarized in Table \ref{table:kappa}.
\begin{table*}[!htbp]
\caption{Performance Metrics of TBM vs POMS under Varying Maintenance Ratios}
\centering
\adjustbox{width=0.9\textwidth}{% NOTE: N = 10 in this table.
\begin{tabular}{ c l C{2.8 cm} C{2.9 cm} C{2.3 cm} C{2.8 cm} C{2.5 cm} C{2 cm} } 
\toprule
    \textbf{\Large $\boldsymbol{\frac{Y^f}{Y^p}}$} &\textbf{Policy}& \textbf{Power Production \scriptsize{($MW$)}} & \textbf{Net Power Transaction \scriptsize{($MW$)}} & \textbf{Availability} & \textbf{Water Spillage \scriptsize{($Hm^3$)} } &\textbf{Operational Revenue}& \textbf{Net Profit} \\ \midrule \midrule
	\multirow{2}{*}{\textbf{1}} & \textbf{TBM}& {12.24 \textbf{M}} & {8.10 \textbf{M}} & {95.87\%} & {32.41} & {\$191.18 \textbf{M}}&{\$183.10 \textbf{M}}\\
	& \textbf{POMS}& {12.30 \textbf{M}} & {8.16 \textbf{M}} & {97.91\%} & {14.11} & {\$194.23 \textbf{M}}&{\$190.37 \textbf{M}} \\
% 	\\ \\\\[-3.6\medskipamount]
\hdashline
% \\\\[-3.3\medskipamount]
	\multirow{2}{*}{\textbf{2}} & \textbf{TBM} & {12.22 \textbf{M}} & {8.08 \textbf{M}} & {95.25\%} & {87.85} & {\$190.23 \textbf{M}} &{\$182.19 \textbf{M}}\\	
	& \textbf{POMS} & {12.30 \textbf{M}} & {8.16 \textbf{M}} & {97.83\%} & {27.01} &  {\$194.06 \textbf{M}}&{\$190.24 \textbf{M}} \\ 
% 	\\\\[-3.6\medskipamount]
	\hdashline
% 	\\\\[-3.3\medskipamount]
	\multirow{2}{*}{\textbf{3}}	& \textbf{TBM} & {12.20 \textbf{M}} & {8.06 \textbf{M}} & {94.55\%} & {155.31} & {\$189.35 \textbf{M}}&{\$181.41 \textbf{M}} \\
	& \textbf{POMS} & {12.29 \textbf{M}} & {8.15 \textbf{M}} & {97.59\%} & {27.01}& {\$193.77 \textbf{M}}&{\$189.75 \textbf{M}} \\ \bottomrule
\end{tabular}}\label{table:kappa}
\end{table*}

Both models exhibit a reduction in availability rates as the duration ratio of corrective maintenance increases. However, POMS sees a minimal reduction of \( \geq 0.33\% \) compared to TBM's \( \geq 1.38\% \) due to lower failure rates of {\ff hydro generators}. The total water spillage of hydropower systems under both models also increases. Specifically, the POMS model faces an additional spillage of 12.9 \(Hm^3\) whereas TBM experiences an additional water spillage of 122.9 \(Hm^3\) when the maintenance ratio increases from 1 to 3. The other impacted performance metrics are the total power production and net power transactions. The power production in both models decreases as a result of higher downtime of {\ff hydro generators}. However, TBM faces a reduction of \( \geq 0.32\% \) in comparison to POMS with only  \( \geq 0.08\% \) power production reduction. The lower power production leads to lower power transaction capability of hydropower systems in the day-ahead market. However the impact on the TBM model with a reduction of  \( \geq 0.49\% \) is harsher in comparison to the POMS model with a reduction of only \( \geq 0.12\% \).
  
It is evident that the POMS model continuously outperforms TBM across all key performance indicators. The sensor-driven nature of POMS becomes increasingly beneficial as corrective maintenance durations increase. Despite the changing \(\frac{Y^f}{Y^p}\) ratio, POMS consistently outperforms TBM in all key performance metrics, including total production (\(\geq 0.50\%\)), net power transaction (\(\geq 0.75\%\)), operational revenue (\(\geq 1.6\%\)), and net profit (\(\geq 3.98\%\)). Unlike TBM, POMS can effectively compensate for longer maintenance downtimes, thereby minimizing negative impacts on power production, revenue, and profit. This underscores the importance of sensor data in optimizing {\ff hydro generators} maintenance schedules, irrespective of the flexibility in corrective maintenance processes.

\subsubsection{\textbf{Experiment III: Transition to Prognostics-driven O\&M Scheduling}}\label{sec:Results_Crew}
The full-scale integration of CM technologies across hydropower systems can be both logistically challenging and financially burdensome \citep{de2016review}. The objective of this experiment is to investigate the benefits of partial deployment of the sensor-driven POMS model for O\&M scheduling of hydropower systems. More specifically, we analyze the performance metrics under hybrid maintenance policies that are mixtures of the TBM and POMS for scheduling O\&M activities of hydropower systems.
     \begin{figure}[!htbp]
    \centering
    \includegraphics[scale=0.5]{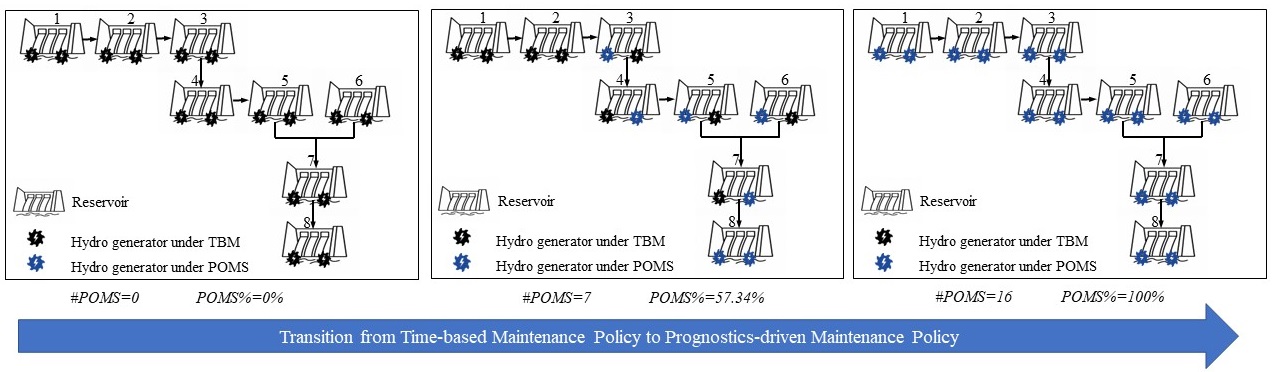}
\caption{Experiment III: Transition from Time-based to prognostics-driven maintenance policy}
\label{fig:Exp3_Example}
\end{figure}  
We quantify the extent of adopting the prognostics-driven O\&M scheduling, i.e., utilizing degradation sensor data of units for failure risk predictions and O\&M scheduling, based on two factors:
\begin{enumerate}[label=(\roman*)]
\item Number of {\ff Hydro Generators} Managed by POMS (\#POMS): This metric captures the number of {\ff hydro generators} equipped with the CM tools and their maintenance schedules are controlled by the POMS model. The range is Min=0, Max=16.
\item Share of Total Power Production Capacity Managed by POMS (\%POMS): This metric captures the percentage of the total power production capacity of hydropower systems managed by POMS. The range is Min=0\%, Max=100\%.
\end{enumerate}
For instance, when \%POMS=0.00\%, it indicates that none of the 16 {\ff hydro generators} are managed using POMS policies, aligning with \#POMS=0. In other words, maintenance actions of the entire fleet of {\ff hydro generators} are scheduled periodically, without utilizing sensor information, using the TBM model. Conversely, \%POMS=100.00\% implies that all {\ff hydro generators} are equipped with sensors and are managed under the POMS model, aligning with \#POMS=16. It is crucial to note two other important factors here:
\begin{enumerate}[]
\item The relationship between the \#POMS and POMS\% is nonlinear and complex. More specifically, an increase in the POMS\% does not necessarily correspond to an increase in the number of {\ff hydro generators} managed by the POMS, i.e., \#POMS. A high POMS\% could be due to a few, high-capacity {\ff hydro generators} coming under POMS management. 
% For example, at 26.33\% capacity under POMS, there are four high-capacity {\ff hydro generators} managed by the POMS model while at 17.52\%, there are ten low-capacity {\ff hydro generators} managed by the POMS model. 
\item The topology of the hydropower systems also impacts the overall performance of the entire system as the functioning of downstream reservoirs is dependent on the operations of those upstream. As a result, the location of {\ff hydro generators} under POMS also plays a significant role in improving the performance metrics.
\end{enumerate}
The interplay between the proportion of {\ff hydro generators} managed under POMS and the spatial and temporal layout of hydropower systems offers a complex but highly relevant avenue for optimization. Strategic POMS implementation could not only depend on the percentage of total production under its aegis but also on where those {\ff hydro generators} are situated in the hydropower systems. The partial deployment of POMS thus has the potential to manifest in different ways.

\begin{table*}[!htbp]
\caption{Impact of POMS adoption on reliability, availability, and maintenance-related metrics}
\centering
\adjustbox{width=0.65\textwidth}{%
\begin{tabular}{C{2cm} C{2cm} C{2.5cm} C{2cm} C{2.5cm} C{2.5cm}}
\toprule
\textbf{\# POMS} & \textbf{\% POMS} & \textbf{Preventive Actions} & \textbf{Failures} & \textbf{Unused Life \scriptsize{(weeks)}} & \textbf{Availability} \\ \midrule \midrule
\textbf{0} & \textbf{0.00\%} & 15 & 6.7 & 59.30 & 95.25\% \\
\textbf{1} & \textbf{8.93\%} & 14.60 & 6.3 & 55.31 & 95.65\% \\
\textbf{2} & \textbf{22.44\%} & 14.60 & 6.5 & 57.73 & 95.31\% \\
\textbf{3} & \textbf{13.09\%} & 13.70 & 6.1 & 54.94 & 95.44\% \\
\textbf{4} & \textbf{26.33\%} & 13.80 & 6 & 53.21 & 95.55\% \\
\textbf{6} & \textbf{35.99\%} & 14.10 & 4.6 & 44.14 & 96.40\% \\
\textbf{7} & \textbf{57.34\%} & 13.70 & 4.5 & 41.94 & 96.59\% \\
\textbf{8} & \textbf{43.56\%} & 14.50 & 3.7 & 42.69 & 96.81\% \\
\textbf{8} & \textbf{75.00\%} & 13.5 & 4.6 & 41.28 & 96.60\% \\
\textbf{9} & \textbf{30.10\%} & 13.40 & 4 & 40.78 & 96.72\% \\
\textbf{9} & \textbf{51.13\%} & 12.80 & 4.3 & 38.61 & 96.49\% \\
\textbf{9} & \textbf{70.59\%} & 12.9 & 4.5 & 38.55 & 96.47\% \\
\textbf{11} & \textbf{66.52\%} & 12.9 & 3.6 & 33.41 & 96.98\% \\
\textbf{13} & \textbf{85.56\%} & 13.9 & 2.5 & 27.50 & 97.38\% \\
\textbf{14} & \textbf{89.94\%} & 13.10 & 2.4 & 26.44 & 97.56\% \\
\textbf{16} & \textbf{100.0\%} & 12.30 & 1.9 & 21.14 & 97.83\% \\ \bottomrule
\end{tabular}}\label{table:Adopt_Maintenance}
\end{table*}

\begin{table*}[!htbp]
\caption{Impact of POMS adoption on operations and expenditures}
\centering
\adjustbox{width=0.8\textwidth}{%
\begin{tabular}{C{2cm} C{2cm} C{3.5cm} C{3.5cm} C{3cm} C{3cm}  C{2cm}}
\toprule
\textbf{\# POMS} & \textbf{\% POMS} & \textbf{Power Production \scriptsize{($MW$)}} & \textbf{Net Power Transaction \scriptsize{($MW$)}} & \textbf{Operational Revenue} & \textbf{Maintenance Cost} & \textbf{Net Profit} \\ \midrule \midrule
\textbf{0} & \textbf{0.00\%} & 12.22 \textbf{M} & 7.97 \textbf{M} & \$190.23 \textbf{M} & \$8.04 \textbf{M} & \$182.19 \textbf{M} \\
\textbf{1} & \textbf{8.93\%} & 12.25 \textbf{M} & 7.95 \textbf{M} & \$191.37 \textbf{M} & \$7.72 \textbf{M} & \$183.65 \textbf{M} \\
\textbf{2} & \textbf{22.44\%} & 12.24 \textbf{M} & 7.38 \textbf{M} & \$190.69 \textbf{M} & \$7.64 \textbf{M} & \$183.05 \textbf{M} \\
\textbf{3} & \textbf{13.09\%} & 12.23 \textbf{M} & 7.88 \textbf{M} & \$190.40 \textbf{M} & \$7.22 \textbf{M} & \$183.18 \textbf{M} \\
\textbf{4} & \textbf{26.33\%} & 12.24 \textbf{M} & 8.07 \textbf{M} & \$190.22 \textbf{M} & \$7.08 \textbf{M} & \$183.14\textbf{M} \\
\textbf{6} & \textbf{35.99\%} & 12.30 \textbf{M} & 8.08 \textbf{M} & \$192.32 \textbf{M} & \$6.02 \textbf{M} & \$186.30 \textbf{M} \\
\textbf{7} & \textbf{57.34\%} & 12.27 \textbf{M} & 8.06 \textbf{M} & \$192.77 \textbf{M} & \$5.86 \textbf{M} & \$186.91 \textbf{M} \\
\textbf{8} & \textbf{43.56\%} & 12.24 \textbf{M} & 7.96 \textbf{M} & \$191.32 \textbf{M} & \$5.78 \textbf{M} & \$185.54 \textbf{M} \\
\textbf{8} & \textbf{75.00\%} & 12.28 \textbf{M} & 8.16 \textbf{M} & \$193.24 \textbf{M} & \$6.06 \textbf{M} & \$187.18 \textbf{M} \\
\textbf{9} & \textbf{30.10\%} & 12.28 \textbf{M} & 8.15 \textbf{M} & \$192.29 \textbf{M} & \$5.64 \textbf{M} & \$186.65 \textbf{M} \\
\textbf{9} & \textbf{51.13\%} & 12.22 \textbf{M} & 8.06 \textbf{M} & \$191.18 \textbf{M} & \$5.84 \textbf{M} & \$185.34 \textbf{M} \\
\textbf{9} & \textbf{70.59\%} & 12.27 \textbf{M} & 8.15 \textbf{M} & \$193.02 \textbf{M} & \$5.78 \textbf{M} & \$187.24 \textbf{M} \\
\textbf{11} & \textbf{66.52\%} & 12.26 \textbf{M} & 8.07 \textbf{M} & \$192.46 \textbf{M} & \$5.30 \textbf{M} & \$187.16 \textbf{M} \\
\textbf{13} & \textbf{85.56\%} & 12.30 \textbf{M} & 8.07 \textbf{M} & \$193.90 \textbf{M} & \$4.62 \textbf{M} & \$189.28 \textbf{M} \\
\textbf{14} & \textbf{89.94\%} & 12.27 \textbf{M} & 8.16 \textbf{M} & \$192.96 \textbf{M} & \$4.30 \textbf{M} & \$188.66 \textbf{M} \\
\textbf{16} & \textbf{100.0\%} & 12.30 \textbf{M} & 8.16 \textbf{M} & \$194.06 \textbf{M} & \$3.82 \textbf{M} & \$190.40 \textbf{M} \\ \bottomrule
\end{tabular}}\label{table:Adopt_Operations}
\end{table*}

\begin{figure*}[htbp!]
    \centering 
    \includegraphics[width=0.45\columnwidth]{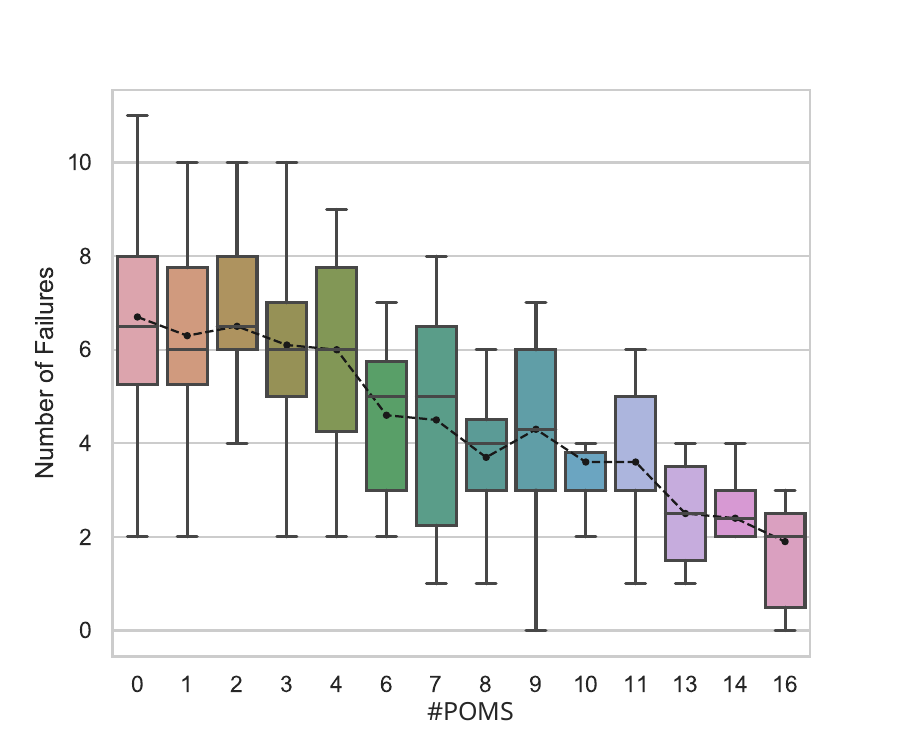}
    % \vspace{-3mm}
    \includegraphics[width=0.45\columnwidth]{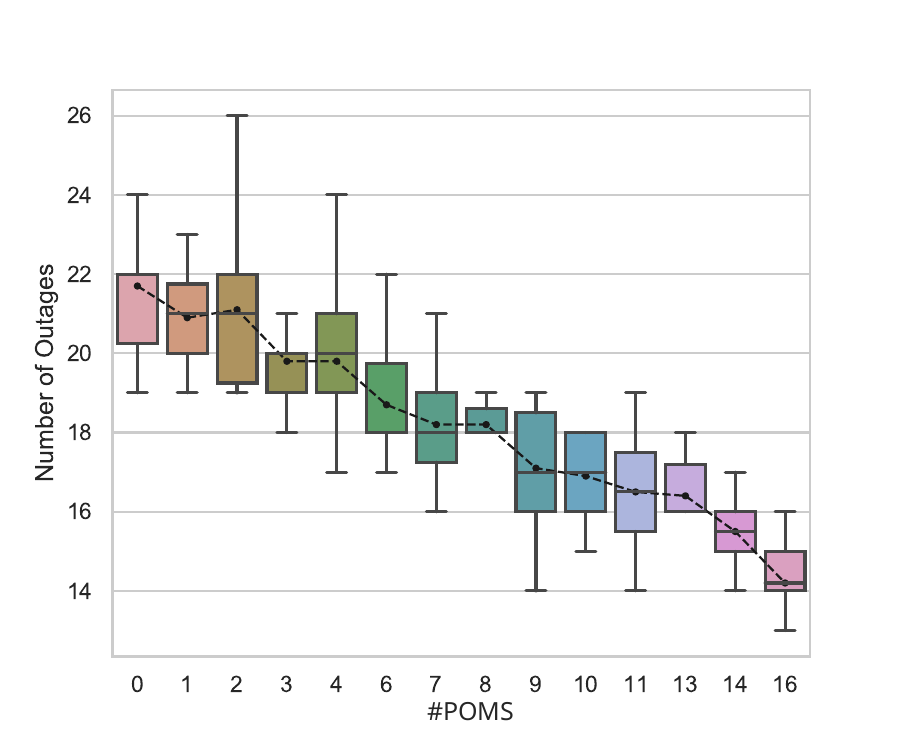}  
    % \vspace*{-1.2cm}
             \caption{Impact of POMS adoption on failures (left) and outages (right) of hydro generators}                      \label{fig:MaintenanceActions} 
        % \vspace{-19mm}
\end{figure*}

\begin{figure*}[htbp!]
    \centering 
    \includegraphics[width=0.45\columnwidth]{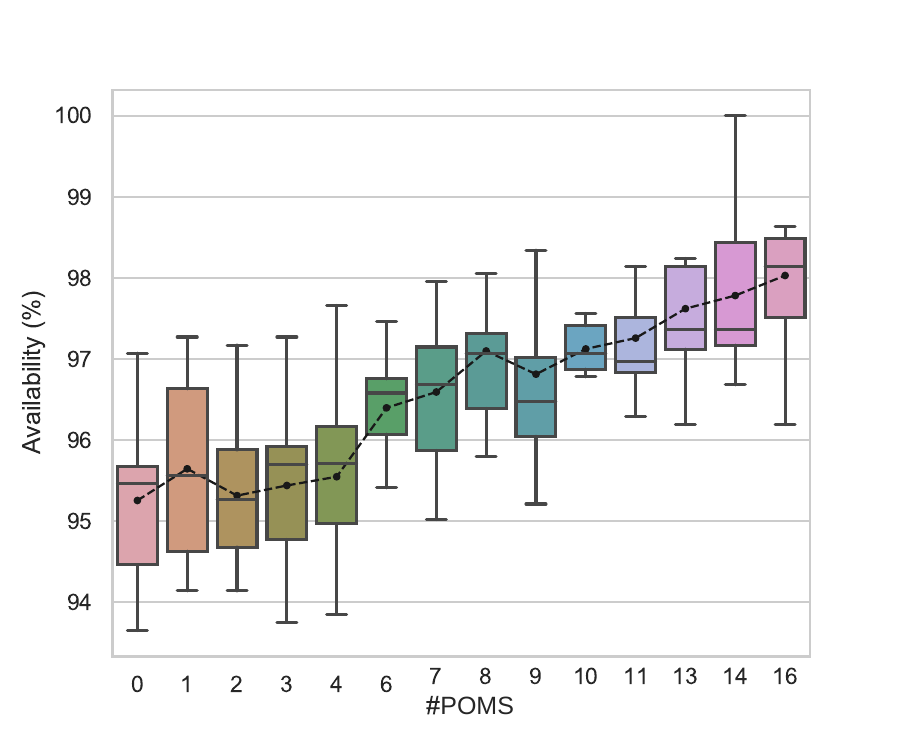}
    % \vspace{-3mm}
    \includegraphics[width=0.45\columnwidth]{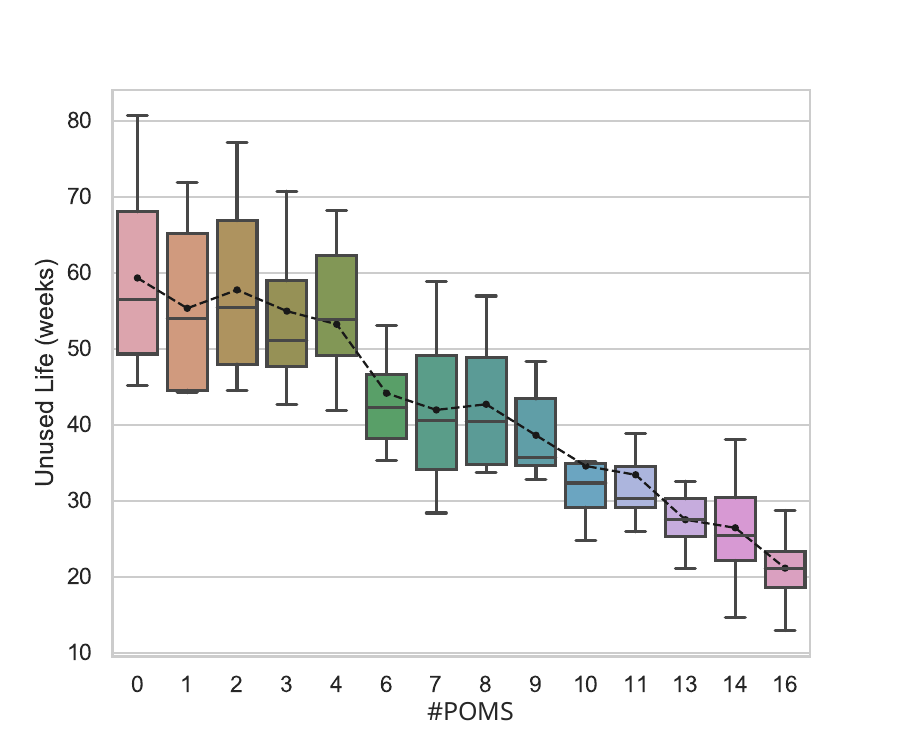}  
    % \vspace*{-1.2cm}
             \caption{Impact of POMS adoption on the availability (left) and unused life (right) of hydro generators} 
                     \label{fig:efficiency} 
        % \vspace{-19mm}
\end{figure*}

\begin{figure*}[htbp!]
    \centering 
    \includegraphics[width=0.45\columnwidth]{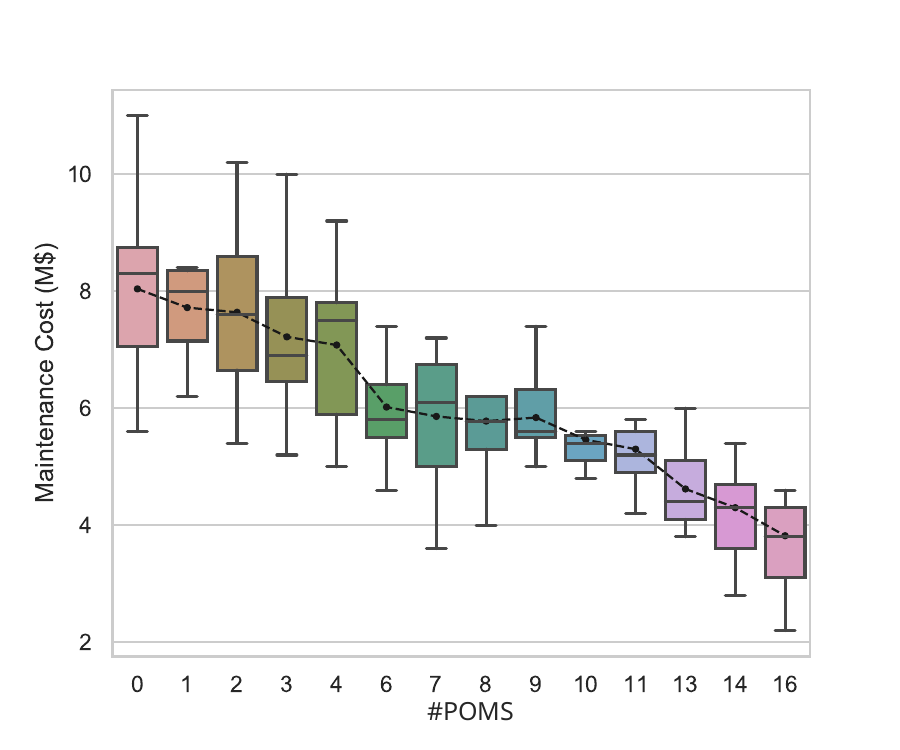}
    % \vspace{-3mm}
    \includegraphics[width=0.45\columnwidth]{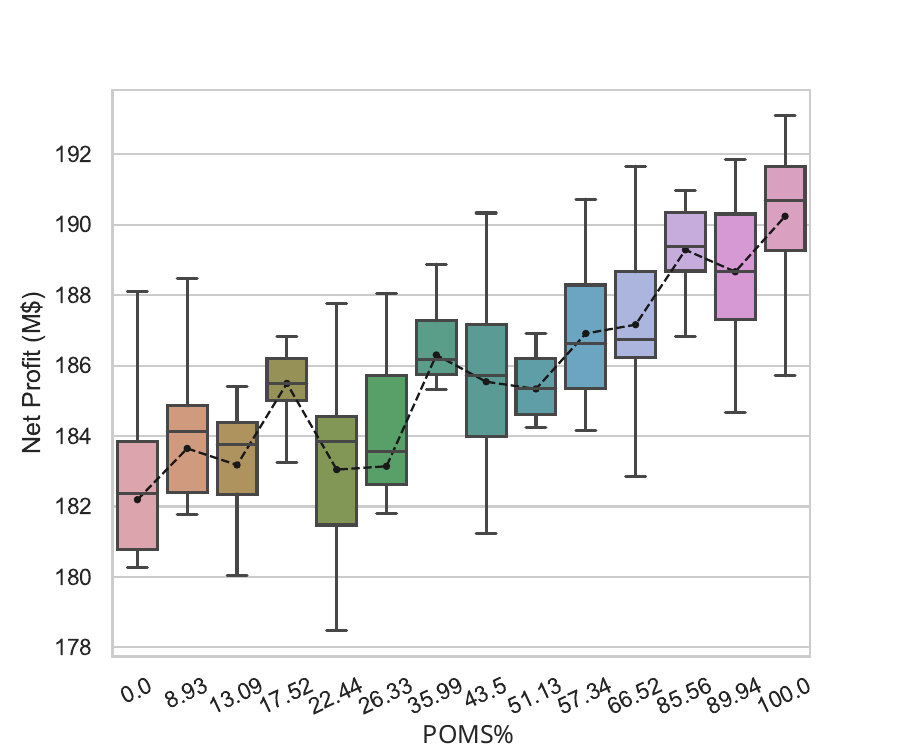}  
    % \vspace*{-1.2cm}
    \caption{Impact of POMS adoption on  maintenance costs (left) and net profit (right) of a hydropower system} 
                     \label{fig:expenditures} 
        % \vspace{-19mm}
\end{figure*}

% To enhance generalizability, 
In this experiment, we gradually increase the number of {\ff hydro generators} under POMS, i.e., \#POMS, from 0 to 16, and randomly choose \#POMS assets to maintain using the POMS policy. We then calculate the proportion of total production capacities, i.e., \%POMS, under POMS policy. The first two columns in tables \ref{table:Adopt_Maintenance} and \ref{table:Adopt_Operations} show the level
of adopting the POMS policy based on these two factors. For each level, we replicate the O\&M scheduling problem with ten different degradation signals and initial ages randomly assigned to {\ff hydro generators} and average the recorded performance metrics. Table \ref{table:Adopt_Maintenance} and \ref{table:Adopt_Operations} display the averaged performance metrics for different levels of POMS adoption for the entire hydropower systems.  We also employ box plots to demonstrate the distribution of the performance metrics across different adoption levels (figures \ref{fig:MaintenanceActions} to \ref{fig:expenditures}). The dashed line in all plots connects the average values. 
% We observe a consistent upward trend across many performance metrics as the \#POMS increases. 
We highlight our observations in the following:

\emph{Enhancing Reliability and Maintenance Metrics:}
 As the \#POMS increases, we see reductions in the average number of preventive maintenance actions and failures of {\ff hydro generators}. Hence, the number of outages of {\ff hydro generators} decreases. The left and right plots in figure \ref{fig:MaintenanceActions} demonstrate the dispersion of the number of failures and outages of {\ff hydro generators} from \#POMS=0 to \#POMS=16. As the number of {\ff hydro generators} managed by POMS policy increases, the average number of failure instances decreases. By increasing the \#POMS, the O\&M scheduling increasingly relies on sensor-driven failure risk predictions of {\ff hydro generators}, enabling it to identify critical conditions in assets and schedule maintenance actions accordingly. This also manifests as a narrowing dispersion in the number of failures and outages, indicating significant improvements in the reliability of O\&M scheduling as the Prognostic-driven POMS policy scales.

\emph{Enhancing Availability and Unused Life Metrics:} Increased reliability of {\ff hydro generators} leads to a higher availability rate of assets as the \#POMS increases.  The left plot in Figure \ref{fig:efficiency} shows an upward improving trend of the availability rates of {\ff hydro generators} as the \#POMS increases. Moreover, increasing the number of assets managed under POMS favors the optimum use of {\ff hydro generators} lifetimes (Figure  \ref{fig:efficiency}, right plot). This improvement stems from the increasing insights from CM and sensor data which enhances the lifetime utilization of assets. In fact, as the number of {\ff hydro generators} equipped with CM tools increases, the model moves from a conservative maintenance scheduling agenda to a liberal maintenance policy approach, enabling more efficient usage of {\ff hydro generators} over the 64-week period. 

\emph{Reducing Expenditures:}
Lastly, in terms of expenditures, a noticeable decrease in maintenance costs is observed as we transition from managing 0 to all 16 {\ff hydro generators} under the POMS framework. This is clearly illustrated in Figure \ref{fig:expenditures} (left plot), where both the average and variability of maintenance costs show improvements. Simultaneously, the proportion of total power production capacity under POMS, denoted as POMS\%, has a direct influence on power production, net power transactions, and as a result profitability. This relationship is highlighted in Figure \ref{fig:expenditures} (right plot), which showcases the average and variability of net profit for a range from POMS\%=0\% to POMS\%=100\%. The results underscore the substantial profitability gains achievable through even partial POMS implementation. The increased visibility on failure risks provided by POMS allows for more agile responses to operational uncertainties and minimizes both the frequency and impact of maintenance-related disruptions. In essence, higher POMS adoption leads to a more reliable maintenance strategy, offering greater scheduling flexibility for {\ff hydro generators} in uncertain operational environments and, ultimately, boosting profitability.
\section{Conclusion}
In this study, we develop a prognostics-driven CBM framework that leverages degradation sensor data and predictive analytics for joint maintenance and operations scheduling in hydropower systems. 
% We utilize Bayesian statistics to continuously derive the most accurate failure risk predictions of {\ff hydro generators} using sensor data. 
We model the joint maintenance and operations scheduling problem of hydropower systems as a two-stage mixed-integer stochastic model and demonstrate the feasibility of integrating the sensor-driven failure risk predictions into the O\&M scheduling problem. To increase the scalability of the framework, we develop a two-level decomposition-based solution algorithm that efficiently solves the O\&M scheduling problem for large-scale cases. We design an experimental framework and evaluate the performance of the prognostics-driven CBM framework and the benchmark time-based periodic maintenance policy under different problem settings. The comprehensive experimental results demonstrate significant benefits of the prognostics-driven framework for the joint O\&M scheduling problem. Specifically, experiments indicate a significant boost in terms of reliability, availability, effective usage of assets lifetime, and expenditures of hydropower systems.  These quantifiable improvements underscore the substantial impact of the developed prognostics-driven O\&M scheduling framework for enhancing the overall performance of hydropower systems. Moreover, experiments prove that even partial deployment of the prognostics-driven framework leads to meaningful gains across various performance metrics of hydropower systems.

Overall, this paper contributes to the advancement of data-driven O\&M scheduling in energy systems and showcases significant benefits that go beyond hydropower system implementations. %, setting the standard for O\&M scheduling across various energy systems.
The proposed framework and solution algorithm can be readily applied to other complex systems with interdependent {\ff generation assets} and uncertain operating conditions. As a case in point, the developed framework can be used as a basis model to develop new approaches for leveraging prognostic predictions to optimize operations, maintenance and spare part logistics in a wide-range of energy and industrial system applications.

%evaluate the CM tools assignments decisions or align inventory and maintenance crew management decisions. 
% Acknowledgments here
%\ACKNOWLEDGMENT{{\fc This material is based upon work supported by the U.S. Department of Energy, Office of Energy Efficiency and Renewable Energy (EERE), specifically the Water Power Technologies Office (WPTO), under seedling project “An Integrated Framework for Condition Monitoring based Asset Management in Hydropower Fleets:  Adaptive and Scalable Prognostics, Operations, and Maintenance}}

% References here (outcomment the appropriate case)

% CASE 1: BiBTeX used to constantly update the references
%   (while the paper is being written).
\bibliographystyle{informs2014} % outcomment this and next line in Case 1
\bibliography{ref} % if more than one, comma separated

\end{document}